\documentclass[12pt]{article}

\usepackage{epsf,epsfig,amsmath,amsfonts}
\usepackage{graphicx}
\usepackage{xcolor}

\usepackage{amsmath}
\usepackage{delarray}

\newtheorem{thm}{Theorem}

\newtheorem{lem}[thm]{Lemma}

\newtheorem{rem}{Remark}
\newtheorem{ass}{Assumption}

\newcommand{\be}{\begin{eqnarray}}
\newcommand{\ee}{\nonumber \end{eqnarray}}
\newcommand{\ben}{\begin{eqnarray}}
\newcommand{\een}{\end{eqnarray}}

\newcommand{\R}{{\mathbb R}}

\newcommand{\F}{{\mathcal F}}
\newcommand{\T}{{\mathcal T}}
\newcommand{\st}{{\mathcal S}}
\newcommand{\W}{{\mathcal W}}
\newcommand{\domH}{{\mathcal D}_H}
\newcommand{\domh}{{\mathcal D}_h}

\newcommand{\Hvf}{{\mathcal H}}
\newcommand{\hvf}{{\mathfrak h}}

\newcommand{\half}{\frac{1}{2}}

\newcommand{\p}{{\mathbb P}}
\newcommand{\tp}{\tilde{\mathbb P}}
\newcommand{\e}{{\mathbb E}}

\usepackage{color}

\begin{document}

\title{\sc Watermark Options}

\author{
{\sc Neofytos Rodosthenous\footnote{School of Mathematical Sciences,
Queen Mary University of London, London E1 4NS, UK,
\texttt{n.rodosthenous@qmul.ac.uk} }
\ and Mihail Zervos\footnote{Department of Mathematics,
London School of Economics, Houghton Street, London WC2A 2AE,
UK, \texttt{mihalis.zervos@gmail.com}}}}

\maketitle

\begin{abstract}\noindent
We consider a new family of derivatives whose payoffs
become strictly positive when the price of their underlying
asset falls relative to its historical maximum.
We derive the solution to the discretionary stopping problems
arising in the context of pricing their perpetual American
versions by means of an explicit construction of their
value functions.
In particular, we fully characterise the free-boundary
functions that provide the optimal stopping times of these
genuinely two-dimensional problems as the unique solutions
to highly non-linear first order ODEs that have the characteristics
of a separatrix.
The asymptotic growth of these free-boundary functions
can take qualitatively different forms depending on
parameter values, which is an interesting new feature.
\\\\
\noindent
{\em MSC2010 subject classification\/}: 49L20, 60G40.
\\\\
\noindent
{\em JEL subject classification\/}: G13, C61.
\\\\
\noindent {\em Key words and phrases\/}: optimal stopping, running
maximum process, variational inequality, two dimensional free-boundary
problem, separatrix.
\end{abstract}

\section{Introduction}

Put options are the most common financial derivatives that can
be used by investors to hedge against asset price falls as well
as by speculators betting on falling prices.
In particular, out of the money put options can yield
strictly positive payoffs only if the price of their underlying asset
falls below a percentage of its initial value.
In a related spirit, equity default swaps (EDSs) pay out if the
price of their underlying asset drops by more than a given
percentage of its initial value (EDSs were introduced by
J.\,P.\,Morgan London  in 2003 and their pricing was studied
by Medova and Smith~\cite{MS06}).
Further derivatives whose payoffs depend on other quantifications
of asset price falls include the European barrier and binary
options studied by Carr~\cite{C06} and Vecer~\cite{V06}, as
well as the perpetual lookback-American options with floating
strike that were studied by Pedersen~\cite{P00} and
Dai~\cite{D01}.

In this paper, we consider a new class of derivatives whose
payoffs depend on asset price falls relative to their underlying
asset's historical maximum price.
Typically, a hedge fund manager's performance fees are
linked with the value of the fund exceeding a ``high
watermark'', which is an earlier maximum.
We have therefore named this new class of derivatives
``watermark'' options.
Deriving formulas for the risk-neutral pricing of their European
type versions is a cumbersome but standard exercise.
On the other hand, the pricing of their American type versions
is a substantially harder problem, as expected.
Here, we derive
the complete solution to the optimal stopping problems
associated with the pricing of their perpetual American
versions.

To fix ideas, we assume that the underlying asset price process
$X$ is modelled by the geometric Brownian motion given by
\ben
dX_t = \mu X_t \, dt + \sigma X_t \, dW_t , \quad X_0 = x > 0 ,
\label{X}
\een
for some constants $\mu$ and $\sigma \neq 0$, where $W$ is
a standard one-dimensional Brownian motion.
Given a point $s \geq x$, we denote by $S$ the running maximum
process defined by
\ben
S_t = \max \left\{ s, \max _{0 \leq u \leq t} X_u \right\} .
\label{S}
\een
In this context, we consider the discretionary stopping problems
whose value functions are defined by
\begin{gather}
v(x,s) = \sup _{\tau \in \T} \e \left[ e^{-r \tau} \left( \frac{ 
S_\tau^b}{X_\tau^a}  - K \right)^+ {\bf 1} _{\{ \tau < \infty \}} 
\right] , \label{v-intro} \\
u(x,s) = \sup _{\tau \in \T} \e \left[ e^{-r \tau} \bigl( 
S_\tau^b - K X_\tau^a \bigr)^+ {\bf 1} _{\{ \tau < \infty \}} 
\right] \label{u-intro}
\end{gather}
and
\ben
\upsilon (x,s) = \sup _{\tau \in \T} \e \left[ e^{-r \tau}
\frac{S_\tau^b}{X_\tau^a}  {\bf 1} _{\{ \tau < \infty \}} 
\right] , \label{vK=0-intro}
\een
for some constants $a, b, r, K > 0$, where $\T$
is the set of all stopping times.
In practice, the inequality $r \geq \mu$ should hold true.
For instance, in a standard risk-neutral valuation context,
$r>0$ should stand for the short interest rate whereas
$r - \mu \geq 0$ should be the underlying asset's dividend
yield rate.
Alternatively, we could view $\mu > 0$ as the short rate
and $r-\mu \geq 0$ as additional discounting to account
for counter-party risk.
Here, we solve the problems we consider without
assuming that $r \geq \mu$ because such an assumption
would not present any simplifications (see also
Remark~\ref{rem:r-mu}).

Watermark options can be used for the same reasons
as the existing options we have discussed above.
In particular, they could be used to hedge against relative
asset price falls as well as to speculate by betting on prices
falling relatively to their historical maximum.
For instance, they could be used by institutions that are
constrained to hold investment-grade assets only and wish
to trade products that have risk characteristics akin to the ones
of speculative-grade assets (see Medova and Smith~\cite{MS06}
for further discussion that is relevant to such an application).
Furthermore, these options can provide useful risk-management
instruments, particularly, when faced with the possibility of an
asset bubble burst.
Indeed, the payoffs of watermark options increase as
the running maximum process $S$ increases and
the price process $X$ decreases.
As a result, the more the asset price increases before dropping
to a given level, the deeper they may be in the money.

Watermark options can also be of interest as hedging
instruments to firms investing in real assets.
To fix ideas, consider a firm that invests in a
project producing a commodity whose price or demand is
modelled by the process $X$.
The firm's future revenue depends on the stochastic evolution
of the economic indicator $X$, which can collapse for reasons
such as extreme changes in the global economic environment
(e.g., see the recent slump in commodity prices) and~/~or
reasons associated with the emergence of disruptive new
technologies (e.g., see the fate of DVDs or firms such as
Blackbury or NOKIA).
In principle, such a firm could diversify risk by going long
to watermark options.

The applications discussed above justify the introduction
of the watermark options as derivative structures.
These options can also provide alternatives to existing
derivatives that can better fit a range of investors' risk
preferences.
For instance, the version associated with (\ref{vK=0-intro})
effectively identifies with the Russian option
(see Remark~\ref{rem:Med-Rus}).
It is straightforward to check that, if $s \geq 1$, then the
price of the option is increasing as the parameter
$b$ increases, ceteris paribus.
In this case, the watermark option is cheaper than the
corresponding Russian option if $b<a$.
On the other hand, increasing values of the strike price $K$
result in ever lower option prices.
We have not attempted any further analysis in this direction
because this involves rather lengthy calculations and is
beyond the scope of this paper.

The parameters $a, b > 0$ and $K \geq 0$ can
be used to fine tune different risk characteristics.
For instance, the choice of the relative value $b/a$ can
reflect the weight assigned to the underlying asset's
historical best performance relative to the asset's current
value.
In particular, it is worth noting that larger (resp., smaller)
values of $b/a$ attenuate (resp., magnify) the payoff's
volatility that is due to changes of the underlying asset's
price.
In the context of the problem with value function given
by (\ref{vK=0-intro}), the choice of $a$, $b$ can be used to factor
in a power utility of the payoff received by the option's
holder.
Indeed, if we set $a = \tilde{a} q$ and $b = \tilde{b} q$,
then $S_\tau^b / X_\tau^a = \bigl( S_\tau^{\tilde{b}} /
X_\tau^{\tilde{a}} \bigr)^q$ is the CRRA utility with
risk-aversion parameter $1-q$ of the payoff
$S_\tau^{\tilde{b}} / X_\tau^{\tilde{a}}$ received by the
option's holder if the option is exercised at time $\tau$.

From a modelling point of view, the standard in the mathematical
finance literature use of a geometric Brownian motion
as an asset price process is an approximation that is largely
justified by its tractability.
In fact, such a process is not an appropriate model
for an asset price that may be traded as a bubble.
In view of the applications we have discussed above, the
pricing of watermark options when the underlying
asset's price process is modelled by diffusions associated
with local volatility models that have been considered in the
context of financial bubbles (e.g., see Cox and
Hobson~\cite{CH05}) presents an interesting problem for
future research.

The Russian options introduced and studied by Shepp and
Shiryaev~\cite{SS93, SS94} are the special cases that arise
if $a=0$ and $b=1$ in (\ref{vK=0-intro}).
In fact, the value function given by (\ref{vK=0-intro}) identifies
with the value function of a Russian option for any $a, b > 0$
(see Remark~\ref{rem:Med-Rus}).
The lookback-American options with floating strike that
were studied by Pedersen~\cite{P00} and Dai~\cite{D01}
are the special cases that arise for the choices $a=b=1$
and $a=b=K=1$ in (\ref{u-intro}), respectively (see also
Remark~\ref{rem:gen-a=b}).
Other closely related problems that have been studied in
the literature include the well-known perpetual American
put options ($a = 1$, $b = 0$ in (\ref{u-intro})), which were
solved by McKean~\cite{McK65}, the lookback-American
options studied by Guo and Shepp~\cite{GS01} and
Pedersen~\cite{P00} ($a=0$, $b=1$ in (\ref{v})), and the
$\pi$-options introduced and  studied by Guo and
Zervos~\cite{GZ10} ($a < 0$ and $b>0$ in (\ref{v-intro})).

Beyond these references, optimal stopping problems
involving a one-dimensional diffusion and its running
maximum (or minimum) include
Jacka~\cite{J91},
Dubins, Shepp and Shiryaev~\cite{DSS93},
Peskir~\cite{P98},
Graversen and Peskir~\cite{GP98},
Dai and Kwok~\cite{DK05, DK06},
Hobson~\cite{H07}, 
Cox, Hobson and Obloj~\cite{CHO08},
Alvarez and Matom{\"a}ki~\cite{AM14},
and references therein.
Furthermore, Peskir~\cite{P14} solves an optimal stopping
problem involving a one-dimensional diffusion, its running
maximum as well as its running minimum.
Optimal stopping problems with an infinite time horizon
involving spectrally negative L\'{e}vy processes and their
running maximum (or minimum) include
Ott~\cite{O13, O14},
Kyprianou and Ott~\cite{KO14},
and references therein.

In Section~\ref{solution1}, we solve the optimal stopping problem
whose value function is given by (\ref{v-intro}) for $a =1$ and
$b \in \mbox{} ]0, \infty[ \mbox{} \setminus \{ 1 \}$.
To this end, we construct an appropriately smooth solution to the
problem's variational inequality that satisfies the so-called
transversality condition, which is a folklore method.
In particular, we fully determine the free-boundary function
separating the ``waiting'' region from the ``stopping'' region as
the unique solution to a first-order ODE that has the characteristics
of a separatrix.
It turns out that this free-boundary function conforms with the
maximality principle introduced by Peskir~\cite{P98}: it is the
maximal solution to the ODE under consideration that does
not intersect the  diagonal part of the state space's boundary.
The asymptotic growth of this free-boundary function is notably
different in each of the cases $1<b$ and $1>b$, which is a
surprising result (see Remark~\ref{as-gr}).

In Section~\ref{solution2}, we use an appropriate
change of probability measure to solve the optimal stopping
problem whose value function is given by (\ref{u-intro})
for $a =1$ and $b \in \mbox{} ]0, \infty[ \mbox{} \setminus
\{ 1 \}$ by reducing it to the problem studied in
Section~\ref{solution1}.
We also outline how the optimal stopping problem defined
by (\ref{X})--(\ref{v-intro}) for $a = b =1$ reduces to
the the problem given by (\ref{X}), (\ref{S}) and (\ref{u-intro})
for $a = b =1$, which is the one arising in the pricing of a
perpetual American lookback with floating strike option that
has been solved by Pedersen~\cite{P00} and Dai~\cite{D01}
(see Remark~\ref{rem:gen-a=b}).
We then explain how a simple re-parametrisation reduces
the apparently more general optimal stopping problems
defined by (\ref{X})--(\ref{u-intro}) for any $a>0$, $b>0$ to
the corresponding cases with $a=1$, $b>0$ (see
Remark~\ref{rem:gen-a,b}).
Finally, we show that the optimal stopping problem
defined by (\ref{X}), (\ref{S}) and (\ref{vK=0-intro})
reduces to the one arising in the context of pricing a
perpetual Russian option that has been solved by
Shepp and Shiryaev~\cite{SS93,SS94} (see
Remark~\ref{rem:Med-Rus}).

\section{The solution to the main optimal stopping problem}
\label{solution1}

We now solve the optimal stopping problem defined by
(\ref{X})--(\ref{v-intro}) for $a =1$ and $b = p \in \mbox{}
]0, \infty[ \mbox{} \setminus \{ 1 \}$, namely, the problem
defined by (\ref{X}), (\ref{S}) and
\ben
v(x,s) = \sup _{\tau \in \T} \e \left[ e^{-r \tau} \left( \frac{ 
S^p_\tau}{X_\tau}  - K \right)^+ {\bf 1} _{\{ \tau < \infty \}} 
\right] . \label{v}
\een
To fix ideas, we assume in what follows that a filtered
probability space $(\Omega, \F, (\F_t), \p)$ satisfying the
usual conditions and carrying a standard one-dimensional
$(\F _t)$-Brownian motion $W$ has been fixed.
We denote by $\T$ the set of all $(\F_t)$-stopping times.

The solution to the optimal stopping problem that we consider
involves the general solution to the ODE
\ben
\half \sigma^2 x^2 f''(x) + \mu x f'(x) - r f(x) = 0 , \label{ODE}
\een
which is given by
\ben
f(x) = A x^n + B x^m , \label{ODE-sol}
\een
for some $A, B \in \R$, where the constants $m < 0 < n$ are
the solutions to the quadratic equation
\ben
\half \sigma^2 k^2 + \left( \mu - \half \sigma^2 \right) k - r = 0 ,
\label{mn-eq}
\een
given by
\ben
m, n = - \frac{\mu - \half \sigma^2}{\sigma^2} \mp
\sqrt{\left( \frac{\mu - \half \sigma^2}{\sigma^2}
\right)^2 + \frac{2 r}{\sigma^2}} . \label{mn}
\een

We make the following assumption.

\begin{ass} \label{A}  {\rm
The constants $p \in \mbox{} ]0, \infty[ \mbox{} \setminus
\{ 1 \}$, $r, K > 0$, $\mu \in \R$ and $\sigma \neq 0$
are such that
\ben
m+1 < 0 \quad \text{and} \quad n+1 - p > 0 . \label{a,b-ass}
\een
} 
\end{ass}

\begin{rem} \label{rem:r-mu} {\rm
We can check that, given any $r>0$, the equivalences
\be
m+1 < 0 \ \Leftrightarrow \ r + \mu > \sigma^2
\quad \text{and} \quad
1 < n \ \Leftrightarrow \ \mu < r
\ee
hold true.
It follows that Assumption~\ref{A} holds true for a range
of parameter values such that $r \geq \mu$, which is
associated with the applications we have discussed in the
introduction, as well as such that $r < \mu$.
} \mbox{}\hfill$\Box$ \end{rem}

We prove the following result in the Appendix.

\begin{lem} \label{lem:v=infty}
Consider the optimal stopping problem defined by
(\ref{X}), (\ref{S}) and (\ref{v}).
If the problem data is such that either $m+1 > 0$
or $n+1-p < 0$, then $v \equiv \infty$.
\end{lem}

We will solve the problem we study in this section
by constructing a classical solution $w$ to the variational
inequality
\ben
\max \left\{ \half \sigma^2 x^2 w_{xx} (x,s) + \mu x w_x
(x,s) - r w(x,s) , \  \left( \frac{s^p}{x} - K \right)^+ - w(x,s)
\right\} = 0 , \label{HJB}
\een
with boundary condition
\ben
w_s (s,s) = 0 , \label{HJB-BC}
\een
that identifies with the value function $v$.
Given such a solution, we denote by $\st$ and $\W$
the so-called stopping and waiting regions, which are
defined by
\begin{align*}
\st & = \left\{ (x,s) \in \R^2 \, \Big| \, \ 0 < x \leq s
\text{ and } w(x,s) = \left( \frac{s^p}{x} - K \right)^+
\right\} \\
\intertext{and}
\W & = \bigl\{ (x,s) \in \R^2 \mid \ 0 < x \leq s \bigr\}
\setminus \st .
\end{align*}
In particular, we will show that the first hitting time
$\tau_\st$ of $\st$, defined by
\ben
\tau_\st = \inf \bigl\{ t \geq 0 \mid \ (X_t, S_t) \in \st \bigr\}
, \label{opt-tau}
\een
is an optimal stopping time.

To construct the required solution to (\ref{HJB})--(\ref{HJB-BC}),
we first note that it is not optimal to stop whenever the
state process $(X,S)$ takes values in the set
\be
\left\{ (x,s) \in \R^2 \, \bigg| \, \ 0 < x \leq s \text{ and }
\frac{s^p}{x} - K \leq 0 \right\} = \left\{ (x,s) \in \R^2
\, \bigg| \, \ s > 0 \text{ and } \frac{s^p}{K} \leq x \leq s
\right\} .
\ee
On the other hand, the equivalences
\begin{align}
\half \sigma^2 x^2 \frac{\partial^2 (s^p x^{-1} - K)}
{\partial x^2} & + \mu x \frac{\partial (s^p x^{-1} - K)}{\partial x}
- r (s^p x^{-1} - K) \leq 0 \nonumber \\
& \Leftrightarrow \quad
\left( \half \sigma^2 - \left( \mu - \half \sigma^2 \right) - r
\right) \frac{s^p}{x} + rK  \leq 0 \nonumber \\
& \Leftrightarrow \quad
x \leq\frac{(n+1) (m+1) s^p}{nmK} \label{HJB-ineq-obs}
\end{align}
imply that the set
\be
\left\{ (x,s) \in \R^2 \, \bigg| \, \ \frac{(n+1) (m+1) s^p}{nmK}
< x \leq s \right\}
\ee
should be a subset of the continuation region $\W$
as well.
Furthermore, since
\be
\left. \frac{\partial}{\partial s} \left( \frac{s^p}{x} - K \right)
\right| _{x=s} = p s^{p-2} > 0 \quad \text{for all } s > 0 ,
\ee
the half-line $\{ (x,s) \in \R^2 \mid \ x = s > 0 \}$, which
is part of the  state space's boundary, should also be
a subset of $\W$ because the boundary condition
(\ref{HJB-BC}) cannot  hold otherwise.

In view of these observations, we look for a strictly
increasing function $H: [0, \infty [ \, \rightarrow \R$
satisfying  
\ben
H(0) = 0 \quad \text{and} \quad 0 < H(s) < [\Gamma s^p]
\wedge s , \quad \text{for } s > 0 , \label{H-reqs}
\een
where 
\begin{equation}
\Gamma = \frac{(n+1) (m+1)}{nmK} \wedge
\frac{1}{K} = \left. \begin{cases}
\frac{(n+1) (m+1)}{nmK} , & \text{if } \mu < \sigma^2 \\
\frac{1}{K} , & \text{if } \sigma^2 \leq \mu \end{cases}
\right\} , \label{Gamma}
\end{equation}
such that
\begin{align}
\st & = \bigl\{ (x,s) \in \R^2 \mid \ s > 0 \text{ and }
0 < x \leq H(s) \bigr\}  \label{st-W-region} \\
\text{and} \quad
\W & = \bigl\{ (x,s) \in \R^2 \mid \ s > 0 \text{ and }
H(s) < x \leq s \bigr\}
\end{align}
(see Figure~1).


To proceed further, we recall the fact that the function $w(\cdot ,
s)$ should satisfy the ODE (\ref{ODE}) in the interior of the waiting
region $\W$.
Since the general solution to (\ref{ODE}) is given by (\ref{ODE-sol}),
we therefore look for functions $A$ and $B$ such that
\be
w(x,s) = A(s) x^n + B (s) x^m , \quad \text{if } \, (x,s) \in \W .
\ee
To determine the free-boundary $H$ and the functions $A$, $B$,
we first note that the boundary condition (\ref{HJB-BC}) requires
that
\ben
\dot{A}(s) s^n + \dot{B}(s) s^m = 0 . \label{W2BC}
\een
In view of the regularity of the optimal stopping problem's
reward function, we expect that the so-called ``principle
of smooth fit'' should hold true.
Accordingly, we require that $w(\cdot, s)$ should be $C^1$
along the free-boundary point $H(s)$, for $s>0$.
This requirement yields the system of equations
\begin{align}
A(s) H^n(s) + B(s) H^m(s) &= s^p H^{-1} (s) - K ,
\nonumber \\
n A(s) H^n(s) + m B(s) H^m(s) &= - s^p H^{-1} (s)
, \nonumber
\end{align}
which is equivalent to
\begin{align}
A(s) & = \frac{-(m+1) s^p H^{-1} (s) + mK}{n-m} H^{-n} (s) ,
\label{W2A} \\
B(s) & = \frac{(n+1) s^p H^{-1} (s) - nK}{n-m} H^{-m}(s) .
\label{W2B}
\end{align}
Differentiating these expressions with respect to $s$ and substituting
the results for $\dot{A}$ and $\dot{B}$ in (\ref{W2BC}), we can see
that $H$ should satisfy the ODE
\begin{gather}
\dot{H}(s) = \Hvf \bigl( H(s), s \bigr) , \label{H-ODE} \\
\intertext{where}
\Hvf \bigl( \bar{H}, s \bigr) = \frac{p s^{p-1} \bar{H}
\left[ (m+1) \left( s / \bar{H} \right)^n - (n+1) \left( s / \bar{H} 
\right)^m \right]} 
{\bigl[ (m+1) (n+1) s^p - nmK \bar{H} \bigr] 
\left[ \left( s / \bar{H} \right)^n - \left( s / \bar{H} \right)^m 
\right]} . \label{calH}
\end{gather}

In view of (\ref{H-reqs}), we need to determine the solution
to (\ref{H-ODE}) in the domain
\ben
\domH = \bigl\{ (\bar{H}, s) \in \R^2 \mid \ s > 0
\text{ and } 0 < \bar{H} < [\Gamma s^p] \wedge s \bigr\}
\label{DH}
\een
that is such that $H(0) = 0$.
To this end, we cannot just solve (\ref{H-ODE}) with the
initial condition $H(0) = 0$ because $\Hvf (0,0)$ is not
well-defined.
Therefore, we need to consider all solutions to (\ref{H-ODE})
in $\domH$ and identify the unique one that coincides
with the actual free-boundary function $H$.
It turns out that this solution is a separatrix.
Indeed, it divides $\domH$ into two open domains
$\domH^u$ and $\domH^l$ such that the solution to
(\ref{H-ODE}) that passes though any point in $\domH^u$
hits the boundary of $\domH$ at some finite $\hat{s}$,
while, the solution to (\ref{H-ODE}) that passes though
any point in $\domH^l$ has asymptotic growth as $s
\rightarrow \infty$ such that the corresponding solution
to the variational inequality (\ref{HJB}) does not
satisfy the transversality condition that is captured by
the limits on the right-hand side of (\ref{TV-cond})
(see also Remark~\ref{TV-cond-rem} and
Figures~2--3).

To identify this separatrix, we fix any $\delta > 0$ and we
consider the solution to (\ref{H-ODE}) that is such that
\ben
H(s_*) = \delta , \quad \text{for some } s _*> s_\dagger
(\delta) , \label{H-IC}
\een
where $s_\dagger (\delta) \geq \delta$ is the intersection
of the half-line $\{ (x,s) \in \R^2 \mid \ x = \delta \text{ and }
s>0 \}$ with the boundary of $\domH$, which is the unique
solution to the equation 
\ben
\bigl[ \Gamma s_\dagger^p (\delta) \bigr] \wedge
s_\dagger (\delta) = \delta . 
\label{s*constraint}
\een 
The following result, which we prove in the Appendix, is
primarily concerned with identifying $s_* > s_\dagger$
such that the solution to (\ref{H-ODE}) that passes through
$(\delta, s_*)$, namely, satisfies (\ref{H-IC}), coincides
with the separatrix.
Using purely analytical techniques, we have not managed
to show that this point $s_*$ is unique, namely, that there
exists a separatrix rather than a funnel.
For this reason, we establish the result for an interval
$[s_\circ, s^\circ]$ of possible values for $s_*$ such that
the corresponding solution to (\ref{H-ODE}) has the
required properties.
The fact that $s_\circ = s^\circ$ follows immediately from
Theorem~\ref{main}, our main result, thanks to the
uniqueness of the optimal stopping problem's value
function.

\begin{lem} \label{H-lem1}
Suppose that the problem data satisfy Assumption~\ref{A}.
Given any $\delta > 0$, there exist points $s_\circ = s_\circ
(\delta)$ and $s^\circ = s^\circ (\delta)$ satisfying
\ben
\delta \leq s_\dagger (\delta) < s_\circ (\delta) \leq s^\circ
(\delta) < \infty , \label{soo-lem}
\een
where $s_\dagger(\delta) \geq \delta$ is the unique solution
to (\ref{s*constraint}), such that the following statements
hold true for each $s_* \in [s_\circ, s^\circ]$:
\smallskip

\noindent
{\rm (I)}
If $p \in \mbox{} ]0,1[$, then the ODE (\ref{H-ODE}) has a
unique solution $H(\cdot) \equiv H(\cdot; s_*) : \mbox{}
]0, \infty[ \mbox{} \rightarrow \domH$ satisfying (\ref{H-IC})
that is a strictly increasing function such that  
\be
\lim_{s \downarrow 0} H(s) = 0 , \quad H(s) < c s^p
\text{ for all } s > 0 \quad \text{and}  \quad
\lim _{s \rightarrow \infty} \frac{H(s)}{s^p} = c ,
\ee
where
$c = \frac{m+1}{mK} \in \mbox{} ]0, \Gamma[$ (see also
Figure~2).
\smallskip

\noindent
{\rm (II)}
If $p>1$, then the ODE (\ref{H-ODE}) has a unique solution
$H(\cdot) \equiv H(\cdot; s_*) : \mbox{} ]0, \infty[ \mbox{}
\rightarrow \domH$ satisfying (\ref{H-IC}) that is a strictly
increasing function such that  
\be
\lim_{s \downarrow 0} H(s) = 0 , \quad H(s) < c s
\text{ for all } s > 0 \quad \text{and}  \quad
\lim _{s \rightarrow \infty} \frac{H(s)}{s} = c ,
\ee
where
$c = \left( \frac{(m+1) (p-n-1)}{(n+1) (p-m-1)} \right)
^{1/(n-m)} \in \mbox{} ]0,1[$ (see also Figure~3).

\smallskip
\noindent
{\rm (III)}
The corresponding functions $A$ and $B$ defined by
(\ref{W2A}) and (\ref{W2B}) are both strictly positive.
\end{lem}

\begin{rem} \label{TV-cond-rem} {\rm
Beyond the results presented in the last lemma,
we can prove the following:
\smallskip

\noindent
({\em a\/})
Given any $s_* \in \mbox{} ]s_\dagger, s_\circ [$, there
exists a point $\hat{s} = \hat{s} (s_*) \in \mbox{} ]0,
\infty[$ and a function $H(\cdot) \equiv H(\cdot; s_*) :
\mbox{} ]0, \hat{s}[ \mbox{} \rightarrow \domH$ that
satisfies the ODE (\ref{H-ODE}) as well as (\ref{H-IC}).
In particular, this function is strictly increasing and
$\lim _{s \uparrow \hat{s}} H(s) = [\Gamma \hat{s}^p]
\wedge \hat{s}$.
\smallskip

\noindent
({\em b\/})
Given any $s_* > s^\circ$, there exists a strictly increasing
function $H(\cdot) \equiv H(\cdot; s_*) : \mbox{} ]0, \infty[
\mbox{} \rightarrow \domH$ that satisfies the ODE (\ref{H-ODE})
as well as (\ref{H-IC}).

\smallskip
\noindent
Any solution to (\ref{H-ODE}) that is as in ({\em a\/}) does
not identify with the actual free-boundary function $H$
because the corresponding solution to the variational
inequality (\ref{HJB}) does not satisfy the boundary condition
(\ref{HJB-BC}).
On the other hand, any solution to (\ref{H-ODE}) that is as in
({\em b\/}) does not identify with the actual free-boundary
function $H$ because we can show that its asymptotic
growth as $s \rightarrow \infty$ is such that the corresponding
solution $w$ to the variational inequality (\ref{HJB}) does
not satisfy (\ref{w-growth}) and the transversality condition,
which is captured by the limits on the right-hand side of
(\ref{TV-cond}), is not satisfied.
To keep the paper at a reasonable length, we do not expand
on any of these issues that are not really needed for
our main results thanks to the uniqueness of the value
function.
} \mbox{}\hfill$\Box$ \end{rem}

\begin{rem} \label{as-gr} {\rm
The asymptotic growth of $H(s)$ as $s \rightarrow \infty$ takes
qualitatively different forms in each of the cases $p<1$ and
$p>1$ (recall that the parameter $p$ stands for the ratio
$b/a$, where the parameters $a$, $b$ are as in the introduction
(see also Remark~\ref{rem:gen-a,b})).
Indeed, if we denote by $H(\cdot ; p)$ the free-boundary
function to indicate its dependence on the parameter
$p$, then
\be
H(s; p) \simeq \left. \begin{cases} cs^p , & \text{if } p<1
\\ cs , & \text{if } p>1 \end{cases}
\right\} \quad \text{as } s \rightarrow \infty ,
\ee
where $c>0$ is the constant appearing in (I) or (II)
of Lemma~\ref{H-lem1}, according to the case.
Furthermore, $c$ is proportional to (resp., independent of)
$K^{-1}$ if $p<1$ (resp., $p>1$).
} \mbox{}\hfill$\Box$ \end{rem}

We now consider a solution $H$ to the ODE (\ref{H-ODE})
that is as in the previous lemma.
In the following result, which we prove in the Appendix, we
show that the function $w$ defined by
\begin{align}
w & (x,s) \nonumber \\
& = \left. \begin{cases}
s^p x^{-1} - K , & \text{if } (x,s) \in \st \\
A(s) x^n + B(s) x^m , & \text{if } (x,s) \in \W
\end{cases} \right\} \nonumber \\
& = \left. \begin{cases}
s^p x^{-1} - K , & \text{if } 0 < x \leq H(s) \\
\frac{-(m+1) s^p H^{-1}(s) + mK}{n-m} \left( \frac{x}{H(s)}
\right)^n +  \frac{(n+1) s^p H^{-1}(s) - nK}{n-m}
\left( \frac{x}{H(s)} \right)^m , & \text{if } H(s) < x \leq s
\end{cases} \right\} \label{w}
\end{align}
is such that
\begin{align}
(x,s) & \mapsto w(x,s) \text{ is } C^2 \text{ outside }
\bigl\{ (x,s) \in \R^2 \mid \ 0 < x \leq s \text{ and } x = H(s)
\bigr\} , \label{w-reg1} \\
x & \mapsto w(x,s) \text{ is } C^1 \text{ at } H(s) \text{ for all }
s  > 0 , \label{w-reg2}
\end{align}
and satisfies (\ref{HJB})--(\ref{HJB-BC}).

\begin{lem} \label{HJB-lem}
Suppose that the problem data satisfy Assumption~\ref{A}.
Also, consider any $s_* \in \bigl[s_\circ (\delta),
s^\circ (\delta) \bigr]$, where $s_\circ (\delta) \leq s^\circ
(\delta)$ are as in Lemma~\ref{H-lem1}, for some $\delta > 0$,
and let $H(\cdot) = H(\cdot; s_*)$ be the corresponding
solution to the ODE (\ref{H-ODE}) that satisfies (\ref{H-IC}).
The function $w$ defined by (\ref{w}) is strictly positive, it
satisfies the variational inequality (\ref{HJB}) outside the set
$\bigl\{ (x,s) \in \R^2 \mid \ s > 0 \text{ and } x = H(s) \bigr\}$
as well as the boundary condition (\ref{HJB-BC}), and is such
that (\ref{w-reg1})--(\ref{w-reg2}) hold true.
Furthermore, given any $s>0$, there exists a constant
$C = C (s) > 0$ such that
\ben
w(x,u) \leq C (1 + u^\gamma) \quad \text{for all } (x,u)
\in \W \text{ such that } u \geq s , \label{w-growth}
\een
where
\be
\gamma = \left. \begin{cases} n (1-p) , & \text{if } p < 1
\\ p-1 , & \text{if } p > 1 \end{cases} \right\} \in \mbox{} ]0,n[ .
\ee
\end{lem}

We can now prove our main result.

\begin{thm} \label{main}
Consider the optimal stopping problem defined by (\ref{X}),
(\ref{S}) and (\ref{v}), and suppose that the problem data
satisfy Assumption~\ref{A}.
The optimal stopping problem's value function $v$
identifies with the solution $w$ to the variational inequality
(\ref{HJB}) with boundary condition (\ref{HJB-BC})
described in Lemma~\ref{HJB-lem}, and the first hitting
time $\tau _\st$ of the stopping region $\st$, which is
defined as in (\ref{opt-tau}), is optimal.
In particular, $s_\circ (\delta) = s^\circ (\delta)$ for all
$\delta > 0$, where $s_\circ \leq s^\circ$ are as in
Lemma~\ref{H-lem1}.
\end{thm}
{\bf Proof.}
Fix any initial condition $(x,s) \in \st \cup \W$.
Using It\^{o}'s formula, the fact that $S$ increases only in
the set $\{ X_t = S_t \}$ and the boundary condition
(\ref{HJB-BC}), we can see that
\begin{align}
e^{-rT} w (X_T,S_T) = \mbox{} & w (x,s) + \int _0^T e^{-rt}
w_s (S_t,S_t) \, dS_t + M_T \nonumber \\
& + \int _0^T e^{-rt} \left[ \half \sigma^2 X_t^2 w_{xx}
(X_t,S_t) + \mu X_t w_x (X_t,S_t) - r w (X_t,S_t) \right] dt
\nonumber \\
= \mbox{} & w (x,s) + M_T \nonumber \\
& + \int _0^T e^{-rt} \left[ \half \sigma^2 X_t^2 w_{xx}
(X_t,S_t) + \mu X_t w_x (X_t,S_t) - r w (X_t,S_t) \right]
dt , \nonumber
\end{align}
where
\be
M_T = \sigma \int _0^T e^{-rt} X_t w_x (X_t,S_t) \, dW_t .
\ee
It follows that
\begin{align}
e^{-rT} \left( \frac{S_T^p}{X_T} - K \right)^+ = \mbox{}
& w (x,s) + e^{-rT} \left[ \left( \frac{S_T^p}{X_T} - K \right)^+
- w (X_T,S_T) \right] + M_T \nonumber \\
& + \int _0^T e^{-rt} \left[ \half \sigma^2 X_t^2 w_{xx}
(X_t,S_t) + \mu X_t w_x (X_t,S_t) - r w (X_t,S_t) \right] dt .
\nonumber
\end{align}

Given a stopping time $\tau \in \T$ and a localising sequence
of bounded stopping times $(\tau_j)$ for the local martingale
$M$, these calculations imply that
\begin{align}
\e \Biggl[ e^{-r (\tau \wedge \tau_j)} &\left(
\frac{S_{\tau \wedge \tau_j}^p}{X_{\tau \wedge \tau_j}}
- K \right)^+ \Biggr] \nonumber \\
= \mbox{} & w (x,s) + \e \left[ e^{-r (\tau \wedge \tau_j)}
\left[ \left( \frac{S_{\tau \wedge \tau_j}^p}
{X_{\tau \wedge \tau_j}} - K \right)^+ -
w (X_{\tau \wedge \tau_j}, S_{\tau \wedge \tau_j}) \right]
\right] \nonumber \\
& + \e \left[ \int _0^{\tau \wedge \tau_j} e^{-rt} \left[ \half
\sigma^2 X_t^2 w_{xx} (X_t,S_t) + \mu X_t w_x (X_t,S_t)
- r w (X_t,S_t) \right] dt \right] . \label{VT00}
\end{align}
In view of the fact that $w$ satisfies the variational inequality
(\ref{HJB}) and Fatou's lemma, we can see that
\be
\e \left[ e^{-r \tau} \left( \frac{ 
S_\tau^p}{X_\tau}  - K \right)^+ {\bf 1} _{\{ \tau < \infty \}} 
\right] \leq \liminf _{j \rightarrow \infty}
\e \Biggl[ e^{-r (\tau \wedge \tau_j)} \left(
\frac{S_{\tau \wedge \tau_j}^p}{X_{\tau \wedge \tau_j}}
- K \right)^+ \Biggr] \leq w (x,s) ,
\ee
and the inequality
\ben
v (x,s) \leq w (x,s) \quad \text{for all } (x,s) \in \st \cup \W
\label{VT1}
\een
follows.

To prove the reverse inequality and establish the
optimality of $\tau_\st$, we note that, given any
constant $T>0$, (\ref{VT00}) with $\tau = \tau_\st \wedge T$
and the definition (\ref{opt-tau}) of $\tau_\st$ imply that
\begin{align}
\e \Biggl[ e^{-r \tau_\st} &\left( \frac{S_{\tau_\st}^p}
{X_{\tau_\st}} - K \right)^+ {\bf 1}
_{\{ \tau_\st \leq \tau_j \wedge T \}} \Biggr]
= w (x,s) - \e \Bigl[ e^{-r (T \wedge \tau_j)}
w (X_{T \wedge \tau_j}, S_{T \wedge \tau_j})
{\bf 1} _{\{ \tau_\st > \tau_j \wedge T \}} \Bigr]
. \nonumber
\end{align}
In view of (\ref{w-growth}), Lemma~\ref{supSgam-lem}
in the Appendix, the fact that $S$ is an increasing process,
and the dominated and monotone convergence theorems,
we can see that
\begin{align}
\e \Biggl[ e^{-r \tau_\st} & \left( \frac{S_{\tau_\st}^p}
{X_{\tau_\st}} - K \right)^+ {\bf 1} _{\{ \tau_\st < \infty \}}
\Biggr] \nonumber \\
& = w (x,s) -
\lim _{T \rightarrow \infty} \lim _{j \rightarrow \infty}
\e \left[ e^{-r (T \wedge \tau_j)}
w (X_{T \wedge \tau_j}, S_{T \wedge \tau_j})
{\bf 1} _{\{ \tau_\st > \tau_j \wedge T \}} \right]
\nonumber \\
& \geq w (x,s) -
\lim _{T \rightarrow \infty} \lim _{j \rightarrow \infty}
\e \left[ e^{-r (T \wedge \tau_j)} C (1 +
S_{T \wedge \tau_j}^\gamma) \right] \nonumber \\
& = w (x,s) - \lim _{T \rightarrow \infty} \e \left[
e^{-r T} C (1 + S_T^\gamma) \right] \nonumber \\
& = w (x,s) . \label{TV-cond}
\end{align}
Combining this result with (\ref{VT1}), we obtain the
identity $v = w$ and the optimality of $\tau_\st$.
Finally, given any $\delta > 0$, the identity $s_\circ (\delta)
= s^\circ (\delta)$ follows from the uniqueness of the value
function $v$.
\mbox{}\hfill$\Box$

\section{Ramifications and connections with the perpetual
American lookback with floating strike and Russian options}
\label{solution2}

We now solve the optimal stopping problem defined by (\ref{X}),
(\ref{S}) and (\ref{u-intro}) for $a =1$ and $b = p \in \mbox{}
]0, \infty[ \mbox{} \setminus \{ 1 \}$, namely, the problem given
by (\ref{X}), (\ref{S}) and
\ben
u(x,s) = \sup _{\tau \in \T} \e \left[ e^{-r \tau}
\bigl( S_\tau^p - K X_\tau \bigr)^+  {\bf 1} _{\{ \tau < \infty \}} 
\right] , \label{u}
\een
by means of an appropriate change of probability measure
that reduces it to the one we solved in Section~\ref{solution1}.
To this end, we denote
\ben
\tilde{\mu} = \mu + \sigma^2 \quad \text{and} \quad
\tilde{r} = r - \mu , \label{tilmr}
\een
and we make the following assumption that mirrors
Assumption~\ref{A}.

\begin{ass} \label{AN}  {\rm
The constants $p \in \mbox{} ]0, \infty[ \mbox{} \setminus
\{ 1 \}$, $r, K > 0$, $\mu \in \R$ and $\sigma \neq 0$ are
such that
\ben
\tilde{m} + 1 > 0 , \quad \tilde{n} + 1 - p > 0 \quad 
\text{and} \quad \tilde{r} > 0 ,
\label{a,b-assN}
\een
where $\tilde{m} < 0 < \tilde{n}$ are the
solutions to the quadratic equation (\ref{mn-eq}),
which are given by (\ref{mn}) with $\tilde{\mu}$
and $\tilde{r}$ defined by (\ref{tilmr}) in place of 
$\mu$ and $r$.
} 
\end{ass}

\begin{thm} \label{mainN}
Consider the optimal stopping problem defined by
(\ref{X}), (\ref{S}) and (\ref{u}) and suppose that the
problem data satisfy Assumption~\ref{AN}.
The problem's value function is given by
\be
u(x,s) = x v(x,s) \quad \text{for all } s>0 \text{ and }
x \in \mbox{} ]0,s]
\ee
and the first hitting time $\tau _\st$ of the stopping
region $\st$ is optimal, where $v$ is the value
function of the optimal stopping problem defined by
(\ref{X}), (\ref{S}) and (\ref{v}), given by Theorem~\ref{main},
and $\st$ is defined by (\ref{st-W-region}) with 
$\tilde{\mu}$, $\tilde{r}$ defined by (\ref{tilmr}) and the 
associated $\tilde{m}$, $\tilde{n}$ in place of 
$\mu$, $r$ and $m$, $n$.
\end{thm}
{\bf Proof.}
We are going to establish this result by means of
an appropriate change of probability measure.
We therefore consider a canonical underlying
probability space because the problem we solve is
over an infinite time horizon.
To this end, we assume that  $\Omega = C(\R_+)$,
the space of continuous functions mapping $\R_+$
into $\R$, and we denote by $W$ the coordinate
process on this space, which is given by
$W_t (\omega) = \omega (t)$.
Also, we denote by $(\F_t)$ the right-continuous
regularisation of the natural filtration of $W$, which
is defined by $\F_t = \bigcap _{\varepsilon > 0}
\sigma \bigl( W_s , \, s \in [0, t+\varepsilon] \bigr)$,
and we set $\F = \bigvee _{t \geq 0} \F_t$.
In particular, we note that the right-continuity of
$(\F_t)$ implies that the first hitting time of any
open or closed set by an $\R^d$-valued continuous
$(\F_t)$-adapted process is an $(\F_t)$-stopping
time (e.g., see Protter~\cite[Theorems~I.3 and~I.4]{P05}).
Furthermore, we denote by $\p$ (resp., $\tp$) the
probability measure on $(\Omega, \F)$
under which, the process $W$ (resp., the process
$\tilde{W}$ defined by $\tilde{W}_t = - \sigma t
+ W_t$) is a standard $(\F_t)$-Brownian motion
starting from $0$.
The measures $\p$ and $\tp$ are locally
equivalent, and their density process is given by
\be
\left. \frac{d\tp}{d\p} \right| _{\F_T} = Z_T , \quad
\text{for } T \geq 0 ,
\ee
where $Z$ is the exponential martingale defined by
\be
Z_T = \exp \left( - \frac{1}{2} \sigma^2 T + \sigma
W_T \right) .
\ee

Given any $(\F_t)$-stopping time $\tau$, we use the
monotone convergence theorem, the fact that
$\e [Z_T \mid \F_\tau] {\bf 1} _{\{ \tau \leq T \}}
= Z_\tau {\bf 1} _{\{ \tau \leq T \}}$ and the tower
property of conditional expectation to calculate
\begin{align}
\e \left[ e^{- r \tau}  \bigl( S_\tau^p - K X_\tau \bigr)^+
{\bf 1} _{\{ \tau < \infty \}}  \right]
& = \lim_{T \rightarrow \infty} \e \left[ e^{- r \tau}
X_\tau \left( \frac{S_\tau^p}{X_\tau} - K \right)^+
{\bf 1}_{\{ \tau \leq T \}} \right] \nonumber \\
& = \lim_{T \rightarrow \infty} x \, \e \left[
Z_T e^{-\tilde{r} \tau} \left( \frac{S_\tau^p}{X_\tau}
- K \right)^+  {\bf 1} _{\{ \tau \leq T \}} \right]
\nonumber \\
& = \lim_{T \rightarrow \infty} x \, \tilde{\e} \left[
e^{- \tilde{r} \tau} \left( \frac{S_\tau^p}{X_\tau}
- K \right)^+  {\bf 1}_{\{ \tau \leq T \}} \right]
\nonumber\\
& = x \, \tilde{\e} \left[ e^{- \tilde{r} \tau}
\left( \frac{S_\tau^p}{X_\tau} - K \right)^+ 
{\bf 1}_{\{ \tau < \infty \}} \right] , \nonumber
\end{align}
and the conclusions of the theorem follow from the
fact that
\be
dX_t = \tilde{\mu} X_t \, dt + \sigma X_t \,
d\tilde{W}_t 
\ee
and Theorem~\ref{main}.
\mbox{}\hfill$\Box$

\begin{rem} \label{rem:gen-a=b} {\rm
Using a change of probability measure argument such as the
one in the proof of the theorem above, we can see that, if
$p=1$, then the value function defined by (\ref{v}) admits the
expression
\be
v(x,s) = x^{-1} \sup _{\tau \in \T} \overline{\e} \left[
e^{-(r + \mu -\sigma^2) \tau} \left( S_\tau - K X_\tau \right)^+
{\bf 1} _{\{ \tau < \infty \}} \right] ,
\ee
where $X$ is given by
\be
dX_t = (\mu - \sigma^2) X_t \, dt + \sigma X_t \,
d\overline{W}_t , \quad X_0 = x > 0 ,
\ee
and expectations are computed under an appropriate
probability measure $\overline{\p}$ under which
$\overline{W}$ is a standard Brownian motion.
This observation reveals that the optimal stopping problem
defined by (\ref{X}), (\ref{S}) and (\ref{v}) for $p=1$ reduces to
the one arising in the pricing of a perpetual American
lookback with floating strike option, which has been
solved by Pedersen~\cite{P00} and Dai~\cite{D01}.
} \mbox{}\hfill$\Box$ \end{rem}

\begin{rem} \label{rem:gen-a,b} {\rm
If $\hat{X}$ is the geometric Brownian motion given by
\ben
d\hat{X}_t = \hat{\mu} \hat{X}_t \, dt + \hat{\sigma} \hat{X}_t
\, dW_t , \quad \hat{X}_0 = \hat{x} > 0 , \label{hatX}
\een
then
\be
d\hat{X}_t^a = \left( \frac{1}{2} \hat{\sigma}^2 a (a-1) +
\hat{\mu} a \right) \hat{X}_t^a \, dt + \hat{\sigma} a
\hat{X}_t^a \, dW_t , \quad \hat{X}_0^a = \hat{x}^a > 0 ,
\ee
and, given any $\hat{s} \geq \hat{x}$,
\ben
\hat{S}_t = \max \left\{ \hat{s}, \max _{0 \leq u \leq t}
\hat{X}_u \right\} = \left( \max \left\{ \hat{s}^a, \max
_{0 \leq u \leq t} \hat{X}_u^a \right\} \right)^{1/a} .
\label{hatS}
\een
In view of these observations, we can see that the solution
to the optimal stopping problem defined by (\ref{hatX}),
(\ref{hatS}) and
\be
\hat{v} (\hat{x}, \hat{s}) = \sup _{\tau \in \T} \e \left[ e^{-r \tau}
\left( \frac{\hat{S}^b_\tau}{\hat{X}^a_\tau}  - K \right)^+ {\bf 1}
_{\{ \tau < \infty \}} \right] ,
\ee
which identifies with the problem (\ref{X})--(\ref{v-intro})
discussed in the introduction, can be immediately derived
from the solution to the problem given by (\ref{X}), (\ref{S})
and (\ref{v}).
Similarly, we can see that the solution to the optimal stopping
problem defined by (\ref{hatX}), (\ref{hatS}) and
\be
\hat{u} (\hat{x}, \hat{s}) = \sup _{\tau \in \T} \e \left[
e^{-r \tau} \bigl( \hat{S}_\tau^b - K \hat{X}_\tau^a
\bigr)^+ {\bf 1} _{\{ \tau < \infty \}} \right] ,
\ee
which identifies with the problem (\ref{X}),  (\ref{S}) and
(\ref{u-intro}) discussed in the introduction, can be
obtained from the solution to the problem given by (\ref{X}),
(\ref{S}) and (\ref{u}).
In particular,
\be
\hat{v} (\hat{x}, \hat{s}) = v (\hat{x}^a, \hat{s}^a) \quad
\text{and} \quad \hat{u} (\hat{x}, \hat{s}) = u (\hat{x}^a,
\hat{s}^a) ,
\ee
for
\be
\mu = \frac{1}{2} \hat{\sigma}^2 a (a-1) + \hat{\mu} a ,
\quad \sigma = \hat{\sigma} a \quad \text{and} \quad
p = \frac{b}{a} .
\ee
Therefore, having restricted attention to the problems
given by (\ref{X}), (\ref{S}) and (\ref{v}) or (\ref{u}) has
not involved any loss of generality.
} \mbox{}\hfill$\Box$ \end{rem}

\begin{rem} \label{rem:Med-Rus} {\rm
Consider the geometric Brownian motion $X$ given by
(\ref{X}) and its running maximum $S$ given by (\ref{S}).
If $\tilde{X}$ is the geometric Brownian motion
defined by
\be
d\tilde{X}_t = \left( \frac{1}{2} \sigma^2 b (b-1) +
\mu b \right) \tilde{X}_t \, dt + \sigma b \tilde{X}_t
\, dW_t , \quad \tilde{X}_0 = \tilde{x} > 0 ,
\ee
and $\tilde{S}$ is its running maximum given by
\be
\tilde{S}_t = \max \left\{ \tilde{s}, \max _{0 \leq u \leq t}
\tilde{X}_u \right\} , \quad \text{for } \tilde{s} \geq \tilde{x} ,
\ee
then
\be
X = \tilde{X}^{1/b} \text{ and } S = \tilde{S}^{1/b}
\quad \text{if } x = \tilde{x}^{1/b} \text{ and }
s = \tilde{s}^{1/b} .
\ee
It follows that the value function $\upsilon$
defined by (\ref{vK=0-intro}) admits the expression
$\upsilon (x,s) = \tilde{\upsilon} (x^b, s^b)$, where
\be
\tilde{\upsilon} (\tilde{x}, \tilde{s}) = \sup _{\tau \in \T}
\e \left[ e^{- r\tau} \frac{\tilde{S}_\tau}{\tilde{X}_\tau
^{a/b}}  {\bf 1} _{\{ \tau < \infty \}} \right] .
\ee
Using a change of probability measure argument such
as the one in the proof of Theorem~\ref{mainN}, we
can see that
\be
\tilde{\upsilon} (\tilde{x}, \tilde{s}) = \tilde{x}^{-a/b}
\sup _{\tau \in \T} \overline{\e} \left[
e^{- \left( r + \mu a - \frac{1}{2} \sigma^2 a (a+1) \right) \tau}
\tilde{S}_\tau  {\bf 1} _{\{ \tau < \infty \}} \right] ,
\ee
where expectations are computed under an appropriate
probability measure $\overline{\p}$ under which the
dynamics of $\tilde{X}$ are given by
\be
d\tilde{X}_t = \left( \frac{1}{2} \sigma^2 b (b-1) +
\mu b - \sigma^2 ab \right) \tilde{X}_t \, dt + \sigma
b \tilde{X}_t \, d\overline{W}_t , \quad \tilde{X}_0
= \tilde{x} > 0 ,
\ee
for a standard Brownian motion $\overline{W}$.
It follows that the optimal stopping problem defined by
(\ref{X}), (\ref{S}) and (\ref{vK=0-intro}) reduces to the
one arising in the context of pricing a perpetual Russian
option, which has been solved by Shepp and
Shiryaev~\cite{SS93, SS94}.
} \mbox{}\hfill$\Box$ \end{rem}

\section*{Appendix: Proof of results in Section~\ref{solution1}}

We need the following result, the proof of which can be found,
e.g., in Merhi and Zervos~\cite[Lemma~1]{MZ}.

\begin{lem} \label{supSgam-lem}
Given any constants $T>0$ and $\zeta \in \mbox{} ]0,n[$,
there exist constants $\varepsilon_ 1 , \varepsilon_2 > 0$
such that
\be
\e \left[ e^{-rT} S_T^\zeta \right] \leq
\frac{\sigma ^2 \zeta ^2 + \varepsilon _2}{\varepsilon _2}
s^\zeta e^{-\varepsilon _1 T} \quad \text{and} \quad
\e \left[ \sup _{T \geq 0} e^{-rT} S_T^\zeta \right] \leq
\frac{\sigma ^2 \zeta ^2 + \varepsilon
_2}{\varepsilon _2} s^\zeta .
\ee
\end{lem}

\noindent
{\bf Proof of Lemma~\ref{lem:v=infty}.}
Suppose first that $m+1 > 0$.
Since $m<0$ is a solution to the quadratic equation
(\ref{mn-eq}), this inequality implies that
\be
\half \sigma^2 (-1)^2 + \left( \mu - \half \sigma^2 \right)
(-1) - r > 0 .
\ee
It follows that
\begin{align}
v(x,s) & \geq \sup _{t \geq 0} \e \left[ e^{-rt} \left(
\frac{S_t^p}{X_t} - K \right)^+ \right]
\geq \sup _{t \geq 0} \e \Big[ e^{-rt} s^p X_t^{-1} \Big]
- e^{-rt}  K \nonumber \\
& = s^p x^{-1} \sup _{t \geq 0} \exp \left( \left[ \half
\sigma^2 - \left( \mu - \half \sigma^2 \right) - r
\right] t \right) - e^{-rt} K = \infty . \nonumber
\end{align}

On the other hand, the inequality $p-1 > n$ and
the fact that $n>0$ is a solution to the quadratic
equation (\ref{mn-eq}) imply that
\be
\half \sigma^2 (p-1)^2 + \left( \mu - \half \sigma^2 \right)
(p-1) - r > 0 .
\ee
In view of this inequality, we can see that
\begin{align}
v(x,s) & \geq \sup _{t \geq 0} \e \left[ e^{-rt} \left(
\frac{S_t^p}{X_t} - K \right)^+ \right]
\geq \sup _{t \geq 0} \e \Big[ e^{-rt} X_t^{p-1} \Big]
- e^{-rt} K \nonumber \\
& = x^{p-1} \sup _{t \geq 0} \exp \left( \left[ \half \sigma^2
(p-1)^2  + \left( \mu - \half \sigma^2 \right) (p-1) - r \right]
t \right) - e^{-rt} K = \infty , \nonumber
\end{align}
and the proof is complete.
\mbox{}\hfill$\Box$
\bigskip

\noindent
{\bf Proof of Lemma~\ref{H-lem1}.}
Throughout the proof, we fix any $\delta > 0$ and we
denote by  $s_\dagger = s_\dagger (\delta)$ the unique
solution to (\ref{s*constraint}).
Combining the assumption $m+1 < 0$ with the observation
that
\ben
[\Gamma s^p] \wedge s \leq \Gamma s^p
\stackrel{(\ref{Gamma})}{\leq} \frac{(m+1)(n+1)}{nmK}
s^p \label{H(s)-ineq}
\een
and  the definition (\ref{DH}) of $\domH$, we can see that
\begin{align}
(m+1)(n+1) s^p - nmK \bar{H} & < (m+1)(n+1) s^p -  nmK
\bigl( [\Gamma s^p] \wedge s \bigr) \nonumber \\
& \leq 0 \qquad \text{for all } (\bar{H}, s) \in \domH ,
\label{anK-ineq}
\end{align}
which implies that the function $\Hvf$ defined by (\ref{calH})
is strictly positive in $\domH$.
Since $\Hvf$ is locally Lipschitz in the open domain $\domH$,
it follows that, given any $s_* > s_\dagger$, there exist points
$\underline{s}_* \in [0, s_*[$ and $\overline{s}_* \in \mbox{}
]s_*, \infty]$, and a unique strictly increasing function $H(\cdot)
= H(\cdot; s_*) : \mbox{} ]\underline{s}_*, \overline{s}_*[
\mbox{} \rightarrow \domH$ that satisfies the ODE (\ref{H-ODE})
with initial condition (\ref{H-IC}) and such that
\be
\lim _{s \downarrow \underline{s}_*} H(s) , \,
\lim _{s \uparrow \overline{s}_*} H(s) \notin \domH
\ee
(see Piccinini, Stampacchia and Vidossich~\cite[Theorems
I.1.4 and~I.1.5]{PSV}).
Furthermore, we note that uniqueness implies that
\ben
s_*^1 < s_*^2 \quad \Leftrightarrow \quad H(s; s_*^1) >
H(s; s_*^2) \text{ for all } s \in \mbox{} \bigl]
\underline{s}_*^1 \vee \underline{s}_*^2, \overline{s}_*^1
\wedge \overline{s}_*^2 \bigr[ . \label{comp1}
\een

Given a point $s_* > s_\dagger$ and the solution
$H(\cdot; s_*)$ to (\ref{H-ODE})--(\ref{H-IC}) discussed
above, we define
\begin{gather}
h(s) = h(s; s_*) = \left. \begin{cases} s^{-p} H(s; s_*) ,
& \text{if } 0 < p < 1 \\ s^{-1} H(s; s_*) , & \text{if } 1
< p < n+1 \end{cases} \right\} , \quad \text{for }
s \in \mbox{} ]\underline{s}_*, \overline{s}_*[ ,
\label{h} \\
s_\circ = \sup \left\{ s_* > s_\dagger \, \Big| \, \ \sup
_{s \in [s_*, \overline{s}_*[} h(s;s_*) \geq c \right\} \vee
s_\dagger \label{so2} \\
\text{and} \quad s^\circ = \inf \left\{ s_* > s_\dagger
\, \Big| \, \ \sup _{s \in [s_*, \overline{s}_*[} h(s;s_*)
< c \right\} , \label{soo}
\end{gather}
with the usual conventions that $\inf \emptyset = \infty$
and $\sup \emptyset = -\infty$, where $c > 0$ is as in
the statement of the lemma, depending on the case.
We also denote
\begin{align}
\domh^1 & = \left\{ (\bar{h}, s) \in \R^2 \mid \ 
s> 0 \text{ and } 0 < \bar{h} < \Gamma \wedge s^{1-p}
\right\} , \nonumber \\
\domh^2 & = \left\{ (\bar{h}, s) \in \R^2 \mid \ 
s> 0 \text{ and } 0 < \bar{h} < [\Gamma s^{p-1}] \wedge
1 \right\} , \nonumber
\end{align}
and we note that
\ben
(\bar{H}, s) \in \domH
\quad \Leftrightarrow \quad
(s^{-p} \bar{H}, s) \in \domh^1
\quad \Leftrightarrow \quad
(s^{-1} \bar{H}, s) \in \domh^2 .
\een
In particular, we can see that these equivalences and
(\ref{anK-ineq}) imply that
\begin{align}
(m+1)(n+1) - nmK \bar{h} < 0 \quad & \text{for all }
(\bar{h}, s) \in \domh^1 , \label{h-ODE-denom1} \\
(m+1)(n+1)s^{p-1} - nmK \bar{h} < 0 \quad & \text{for all }
(\bar{h}, s) \in \domh^2 , \label{h-ODE-denom2}
\end{align}
while, (\ref{comp1}) implies the equivalence
\ben
s_*^1 < s_*^2 \quad \Leftrightarrow \quad h(s; s_*^1) >
h(s; s_*^2) \text{ for all } s \in \mbox{} \bigl]
\underline{s}_*^1 \vee \underline{s}_*^2, \overline{s}_*^1
\wedge \overline{s}_*^2 \bigr[ . \label{comp2}
\een

In view of definitions (\ref{h})--(\ref{soo}), the required
claims will follow if we prove that
\begin{gather}
\underline{s}_* = 0 \quad \text{and} \quad \overline{s}_*
= \infty \quad \text{for all } s_* \in [s_\circ, s^\circ] ,
\label{us*} \\
s_\dagger < s_\circ \leq s^\circ < \infty , \label{s_^o} \\
h(s; s_\circ) < c \quad \text{for all } s > 0 ,
\label{s-circs-(II).1} \\
\intertext{as well as} 
\limsup_{s \downarrow 0} h(s;s_*) < \infty
\quad \text{and} \quad
\lim _{s \rightarrow \infty} h(s;s_*) = c \quad 
\text{for all } s_* \in [s_\circ, s^\circ] .
\label{s-circs-(II).2}
\end{gather}
To prove that these results are true, we need to differentiate
between the two cases of the lemma: although the main
ideas are the same the calculations involved are remarkably
different (compare Figures~4 and~5).

\underline{{\em Proof of (I)} $(p < 1)$\/.}
In this case, which is illustrated by Figure~4, we calculate
\ben
\dot{h}(s) \equiv \dot{h} (s; s_*) = \hvf \bigl( h(s), s \bigr) 
\quad \text{and} \quad
h(s_*) = \delta s_*^{-p} , \label{h-ODE1}
\een
where
\begin{align}
\hvf (\bar{h}, s) = \mbox{} & - \frac{\bar{h}}{s} \frac{p
\bigl[ n(m+1) - nmK \bar{h} \bigr] \left( s^{1-p} \bar{h}^{-1}
\right)^{n-m}} {\bigl[ (s^{1-p} \bar{h}^{-1})^{n-m}
- 1 \bigr] \bigl[ (m+1)(n+1) - nmK \bar{h} \bigr]} \nonumber \\ 
& + \frac{\bar{h}}{s} \frac{pm(n+1) - pnmK \bar{h}}
{\bigl[ (s^{1-p} \bar{h}^{-1})^{n-m} - 1 \bigr] \bigl[ (m+1)(n+1)
- nmK \bar{h} \bigr]} . \nonumber
\end{align}
In the arguments that we develop, the inequalities
\be
m+1 < 0 < n, \quad 0< p < 1 \quad \text{and} \quad c =
\frac{m+1}{mK} \in \mbox{} ]0, \Gamma[ ,
\ee
which are relevant to the case we now consider,
are worth keeping in mind.
In view of the inequalities
\begin{gather}
n(m+1) - nmK \bar{h} = nmK (c - \bar{h}) \left. \begin{cases}
< 0 , & \text{if } \bar{h} \in \mbox{} ]0, c[ \\ > 0 , & \text{if }
\bar{h} \in \mbox{} ]c, \Gamma[ \end{cases} \right\} ,
\label{h-ODE-denom11} \\
m(n+1) - nmK \bar{h} < (m+1)(n+1) - nmK \bar{h}
\stackrel{(\ref{h-ODE-denom1})}{<} 0 \quad \text{for all }
(\bar{h}, s) \in \domh^1 , \nonumber
\end{gather}
we can see that
\ben
\bigl\{ (\bar{h}, s) \in \domh^1 \mid \ \hvf (\bar{h}, s)
< 0 \bigr\} = \bigl\{ (\bar{h}, s) \in \domh^1 \mid \
\bar{h} < c \text{ and } s > \mathfrak{s} (\bar{h}) \bigr\}
, \label{D1s}
\een
where the function $\mathfrak{s}$ is defined by
\be
\mathfrak{s} (\bar{h}) = \left( \left[ \frac{m(n+1)
- nmK \bar{h}} {n(m+1) - nmK \bar{h}} \right]
^{1/(n-m)} \bar{h}  \right)^{1/(1-p)} , \quad \text{for }
\bar{h} \in \mbox{} ]0,c[ . 
\ee
Furthermore, we calculate
\be
\dot{\mathfrak{s}} (\bar{h}) > 0 \text{ for all } \bar{h}
\in \mbox{} ]0,c[ , \quad \lim _{\bar{h} \downarrow 0}
\mathfrak{s} (\bar{h}) = 0  \quad \text{and} \quad 
\lim _{\bar{h} \uparrow c} \mathfrak{s} (\bar{h}) = \infty ,
\ee
and we note that
\ben
\lim_{\bar{h} \uparrow s^{1-p}} \hvf (\bar{h}, s)
= \infty \quad \text{for all } s \leq c^{1/(1-p)} .
\label{h(h,s)-oo1}
\een
In particular, this limit and (\ref{D1s}) imply that
\be
\mathfrak{s}^{\{ -1 \}} (s) < c \wedge s^{1-p} \leq
\Gamma \wedge s^{1-p} \quad \text{for all } s > 0 ,
\ee
where $\mathfrak{s}^{\{ -1 \}}$ is the inverse function of
$\mathfrak{s}$.
In view of these observations, we can see that
\ben
\hvf (\bar{h}, s) \left. \begin{cases} < 0 &
\text{for all } s>0 \text{ and } \bar{h} \in \mbox{}
]0, \mathfrak{s}^{\{ -1 \}} (s)[ \\ > 0 & \text{for all } s>0
\text{ and } \bar{h} \in \mbox{} \bigl] \mathfrak{s}^{\{ -1 \}}
(s), \Gamma \wedge s^{1-p} \bigr[ \end{cases} \right\} ,
\label{conc11}
\een
as well as that
\ben
\delta s_*^{-p} - \mathfrak{s}^{\{ -1 \}} (s) \left.
\begin{cases} > 0 & \text{for all } s \in [s_\dagger,
s^\dagger[ \\ < 0 & \text{for all } s \in \mbox{}
]s^\dagger, \infty[ \end{cases} \right\} , \label{conc12}
\een
for a unique $s^\dagger = s^\dagger (\delta)
> s_\dagger$.

The conclusions (\ref{conc11})--(\ref{conc12}) imply
immediately that
\be
\sup _{s \in \mbox{} ]s_*, \overline{s}_*[} h(s; s_*) <
\delta s_*^{-p} < c \quad \text{for all } s_* \geq s^\dagger .
\ee
Combining this inequality with (\ref{comp2}) and (\ref{conc11}),
we obtain
\be
s^\circ < s^\dagger \quad \text{and} \quad
\forall s_* \in \mbox{} ]s^\circ, s^\dagger[ , \
\exists! s_{\mathrm m} = s_{\mathrm m} (s_*) : \
h(\cdot; s_*) \left. \begin{cases} \text{is strictly
increasing in } \mbox{} ]s_*, s_{\mathrm m}[ \\
\text{is strictly decreasing in } \mbox{} ]s_{\mathrm m},
\overline{s}_*[ \end{cases} \hspace{-4mm} \right\} .
\ee
In view of this observation, (\ref{comp2}),
(\ref{h(h,s)-oo1}), (\ref{conc11}) and a straightforward
contradiction argument, we can see that
\begin{gather}
s_\circ \in \mbox{} ]s_\dagger, s^\circ] , \quad
h(\cdot; s^\circ) \text{ is strictly increasing} ,
\nonumber \\
\lim _{s \rightarrow \infty} h(s; s_*) = c \text{ for all }
s_* \in [s_\circ, s^\circ] \quad \text{and} \quad
\lim _{s \downarrow 0} h(s;s_*) = 0 \text{ for all } s_*
\in \mbox{} ]s_\dagger, s^\circ] . \nonumber
\end{gather}
It follows that (\ref{us*})--(\ref{s-circs-(II).2}) are
all true.

\underline{{\em Proof of (II)} ($1 < p$)\/.}
In this case, which is illustrated by Figure~5, we calculate
\ben
\dot{h}(s) \equiv \dot{h} (s; s_*) = \hvf \bigl( h(s), s \bigr) 
\quad \text{and} \quad
h(s_*) = \delta s_*^{-1} , \label{h-ODE2}
\een
where
\begin{align}
\hvf (\bar{h}, s) = \mbox{} & \frac{\bar{h}}{s} \frac{s^{p-1} 
\bigl[ (m+1) (p-1-n) - (n+1) (p-1-m) \bar{h}^{n-m} \bigr]}
{\left( 1 - \bar{h}^{n-m} \right)
\left[ (m+1)(n+1) s^{p-1} - nmK \bar{h} \right]} \nonumber \\ 
& + \frac{\bar{h}}{s} \frac{nmK \bar{h}}{(m+1)(n+1)
s^{p-1} - nmK \bar{h}} . \nonumber
\end{align}
In what follows, we use the inequalities
\begin{gather}
m+1 < 0 < n , \quad
1 < p < n+1 \quad \text{and} \quad c = \left(
\frac{(m+1) (p-1-n)}{(n+1) (p-1-m)} \right)^{1/(n-m)}
\in \mbox{} ]0, 1[ \nonumber
\end{gather}
that are relevant to the case we now consider.
In view of (\ref{h-ODE-denom2}) and the inequalities
\begin{align}
(m+1) (p-1-n) & - (n+1)(p-1-m) \bar{h}^{n-m}
\nonumber \\
& = (n+1)(p-1-m) (c^{n-m} - \bar{h}^{n-m})
\left. \begin{cases} > 0 , & \text{if } \bar{h} \in \mbox{}
]0, c[ \\ < 0 , & \text{if } \bar{h} \in \mbox{} ]c, 1[
\end{cases} \right\} ,
\end{align}
we can see that
\ben
\bigl\{ (\bar{h}, s) \in \domh^2 \mid \ \hvf (\bar{h}, s) < 0
\bigr\} = \bigl\{ (\bar{h}, s) \in \domh^2 \mid \
\bar{h} < c \text{ and } s > \mathfrak{s} (\bar{h}) \bigr\} ,
\label{D2s}
\een
where the function $\mathfrak{s}$ is defined by
\be
\mathfrak{s} (\bar{h}) = \left( \frac{- nmK \bigl( 1 -
\bar{h}^{n-m} \bigr)} {(m+1) (p-1-n) - (n+1) (p-1-m)
\bar{h}^{n-m}} \bar{h} \right)^{1/(p-1)} , \quad
\text{for } \bar{h} \in \mbox{} ]0, c[ .
\ee
It is straightforward to check that
\ben
\dot{\mathfrak{s}} (\bar{h}) > 0 \text{ for all } \bar{h} > 0 ,
\quad \lim _{\bar{h} \downarrow 0} \mathfrak{s}
(\bar{h}) = 0  \quad \text{and} \quad 
\lim _{\bar{h} \uparrow c} \mathfrak{s} (\bar{h}) = \infty .
\label{sprop2}
\een

To proceed further, we define
\ben
\Gamma_1 = \frac{(n+1) (m+1)}{nmK}
\quad \text{and} \quad \Gamma_2 = \frac{1}{K} ,
\label{Gamma12}
\een
we note that $\Gamma = \Gamma_1 \wedge
\Gamma_2$, and we observe that
\ben
\lim_{\bar{h} \uparrow \Gamma_1 s^{p-1}} \hvf
(\bar{h}, s) = \infty \quad \text{for all } s \leq \left(
\frac{c}{\Gamma_1} \right)^{1/(p-1)} .
\label{h(h,s)-oo21}
\een
If $\Gamma = \Gamma_2$, then we can use the inequality
$(m+1)(n+1) - nm < 0$, which holds true in this case
(see also (\ref{Gamma})), to verify that
\be
\frac{(m+1)(p-1-n) - (n+1)(p-1-n) \bar{h}^{n-m} + nm
(1-\bar{h}^{n-m})}{\bigl[ (m+1)(n+1) - nm \bigr]
(1-\bar{h}^{n-m})} > p-1 \quad \text{for all }
\bar{h} \in \mbox{} ]0,1[ .
\ee
Using this inequality, we can see that
\ben
\text{if } \Gamma = \Gamma_2 , \text{ then} \quad
\hvf \bigl( \Gamma s^{p-1}, s \bigr) > \frac{d}{ds}
\bigl( \Gamma s^{p-1} \bigr) > 0 \quad \text{for all }
s \leq \left( \frac{c}{\Gamma} \right)^{1/(p-1)} .
\label{h(h,s)-oo22}
\een
Combining (\ref{D2s}) with
(\ref{h(h,s)-oo21})--(\ref{h(h,s)-oo22}), we obtain
\be
\mathfrak{s}^{\{ -1 \}} (s) < \bigl[ \Gamma s^{p-1} \bigr]
\wedge c \quad \text{for all } s > 0 ,
\ee
where $\mathfrak{s}^{\{ -1 \}}$ is the inverse function of
$\mathfrak{s}$.
It follows that
\ben
\hvf (\bar{h}, s) \left. \begin{cases} < 0 &
\text{for all } s>0 \text{ and } \bar{h} \in \mbox{}
]0, \mathfrak{s}^{-1} (s)[ \\ > 0 & \text{for all } s>0
\text{ and } \bar{h} \in \mbox{} \bigl] \mathfrak{s}^{-1}
(s), \Gamma_1 \wedge s^{p-1} \bigr[ \mbox{}
\supseteq \mbox{} \bigl] \mathfrak{s}^{-1} (s), \Gamma
\wedge s^{p-1} \bigr[ \end{cases} \right\} ,
\label{conc21}
\een
as well as that
\ben
\delta s_*^{-1} - \mathfrak{s}^{-1} (s) \left.
\begin{cases} > 0 & \text{for all } s \in [s_\dagger,
s^\dagger[ \\ < 0 & \text{for all } s \in \mbox{}
]s^\dagger, \infty[ \end{cases} \right\} , \label{conc22}
\een
for a unique $s^\dagger = s^\dagger (\delta)
> s_\dagger$.

Arguing in exactly the same way as in Case~(I) above
using (\ref{comp2}) and (\ref{conc21})--(\ref{conc22}),
we can see that (\ref{s_^o})--(\ref{s-circs-(II).2}) are
all true.
Furthermore, (\ref{us*}) follows from
(\ref{sprop2})--(\ref{h(h,s)-oo22}).

\underline{{\em Proof of (III)}\/.}
In view of (\ref{DH}), (\ref{H(s)-ineq}) and
the observation that
\be
m+1 < 0 \quad \Rightarrow \quad
\frac{m+1}{m} \in \mbox{} ]0,1[ ,
\ee
we can see that
\be
H^{-1} (s) > \frac{nK}{n+1} s^{-p} \quad \Leftrightarrow
\quad B(s) > 0 \quad \text{for all } s>0 .
\ee

If $p<1$, then we can use the fact that
\be
s^p H^{-1} (s) > c^{-1} = \frac{mK}{m+1} \quad
\text{for all } s> 0 ,
\ee
which we have established in part~(I) of the lemma,
to see that
\be
A(s) > 0 \quad \Leftrightarrow \quad s^p H^{-1} (s) >
\frac{mK}{m+1} \quad \text{for all } s>0
\ee
is indeed true.

On the other hand, if $p>1$, then we can verify that
\be
\frac{(m+1) \left( s / \bar{H} \right)^n - (n+1) \left(
s / \bar{H} \right)^m} {(m+1) \left[ \left( s / \bar{H}
\right)^n - \left( s / \bar{H} \right)^m \right]} > 1
\quad \text{for all } s > 0 \text{ and } \bar{H} \in
\mbox{} ]0,s[ .
\ee
Using this inequality, we obtain
\be
\Hvf \left( \frac{m+1}{mK} s^p , s \right) > \frac{d}{ds}
\left( \frac{m+1}{mK} s^p \right) \quad \text{for all } s > 0 .
\ee
Combining this calculation with (\ref{DH}),
we can see that
\be
H(s) < \frac{m+1}{mK} s^p \quad \Leftrightarrow \quad
A(s) > 0 \quad \text{for all } s > 0
\ee
because, otherwise, $H$ would exit the domain
$\domH$.
\mbox{}\hfill$\Box$
\bigskip

\noindent
{\bf Proof of Lemma~\ref{HJB-lem}.}
We first note that the strict positivity of $w$ follows
immediately from its definition in (\ref{w}) and
Lemma~\ref{H-lem1}.(III).
To establish (\ref{w-growth}), we fix any $s>0$ and
we note that Lemma~\ref{H-lem1} implies that there
exists a point $\overline{s} \geq s$ such that
\be
H(u) \geq \left. \begin{cases} \frac{1}{2} c u^p ,
& \text{if } p < 1 \\ \frac{1}{2} c u , & \text{if } p > 1
\end{cases} \right\} \quad \text{for all } u \geq
\overline{s} .
\ee
Combining this observation with the fact that $H$ is
continuous, we can see that there exists a constant
$C_1 (s) > 0$ such that
\begin{align}
u^p x^n H^{-(n+1)} (u) & \leq u^{p+n} H^{-(n+1)} (u)
\nonumber \\
& \leq \max _{u \in [s, \overline{s}]} \bigl\{ u^{p+n}
H^{-(n+1)} (u) \bigr\} {\bf 1} _{[s, \overline{s}]} (u)
+ u^{p+n} H^{-(n+1)} (u) {\bf 1} _{[\overline{s}, \infty[}
(u) \nonumber \\
& \leq \left. \begin{cases} C_1 (s) \bigl(1 + u^{n (1-p)}
\bigr) , & \text{if } p < 1 \\ C_1 (s) \bigl(1 + u^{p-1} \bigr)
, & \text{if } p > 1 \end{cases} \right\} \quad \text{for all }
u \geq s \text{ and } x \leq u , \nonumber
\end{align}
and
\begin{align}
u^p x^m H^{-(m+1)} (u) & \leq u^p H^{-1} (u)
\nonumber \\
& \leq \left. \begin{cases} C_1 (s) , & \text{if } p < 1 \\
C_1 (s) \bigl(1 + u^{p-1} \bigr) , & \text{if } p > 1
\end{cases} \right\} \quad \text{for all } u \geq s
\text{ and } x \in [H(u), u] . \nonumber
\end{align}
In view of these calculations, the definition (\ref{w})
of $w$, and the inequalities $m+1 < 0 < n$,
we can see that there exists a constant
$C = C(s) > 0$ such that (\ref{w-growth}) holds true.
For future reference, we note that the second
of the estimates above implies that, given any $s>0$,
\ben
u^{p-1} \leq u^p x^{-1} \leq u^p H^{-1} (u) \leq C_1
(s) (1 + u^{p-1}) \quad \text{for all }
u \geq s \text{ and } x \in [H(u), u] . \label{s/x-est}
\een

By construction, we will prove that the positive function
$w$ is a solution to the variational inequality (\ref{HJB})
with boundary condition (\ref{HJB-BC}) that satisfies
(\ref{w-reg1})--(\ref{w-reg2}) if we show that
\begin{align}
f(x,s) & := \half \sigma^2 x^2 \frac{\partial ^2}{\partial x^2}
\left( \frac{s^p}{x} - K \right) + \mu x \frac{\partial}{\partial x}
\left( \frac{s^p}{x} - K \right) - r \left( \frac{s^p}{x} - K \right)
\nonumber \\
& = \left[ \half \sigma^2 (-1)^2 + \left( \mu - \half \sigma^2
\right) (-1) - r \right] \frac{s^p}{x} + rK
\leq 0 \qquad \text{for all } (x,s) \in \st \label{HJB-S}
\end{align}
and
\ben
g(x,s) := w(x,s) - \frac{s^p}{x} + K \geq 0 \quad \text{for all }
(x,s) \in \W . \label{HJB-W}
\een

\underline{\em Proof of (\ref{HJB-S})\/.}
In view of the assumption that $m+1<0$ and the fact that
$m<0<n$ are the solutions to the quadratic equation
(\ref{mn-eq}) given by (\ref{mn}), we can see that 
\ben
0 > \half \sigma^2 (-1)^2 + \left( \mu - \half \sigma^2 \right)
(-1) - r = \half \sigma^2 (n+1)(m+1) = - \frac{r}{nm}
(n+1)(m+1) . \label{a-quad}
\een
Combining this inequality and identities with the fact that
$H(s) < \Gamma s^p \leq \frac{(n+1) (m+1)}{nmK} s^p$
for all $s > 0$ (see (\ref{H-reqs})--(\ref{Gamma})), we
calculate
\ben
f_x (x,s) = - \left[ \half \sigma^2 (-1)^2 + \left( \mu - \half
\sigma^2 \right) (-1) - r \right] \frac{s^p}{x^2} > 0 \quad
\text{for all } s > 0 \text{ and } x \in \mbox{} ]0, s[ ,
\label{f-deriv}
\een
and
\begin{align}
f \bigl( H(s), s \bigr) & = \left[ \half \sigma^2 (-1)^2 + \left(
\mu - \half \sigma^2 \right) (-1) - r \right] \frac{s^p}{H(s)}
+ rK \nonumber \\
& < \left[ \half \sigma^2 (-1)^2 + \left( \mu - \half
\sigma^2 \right) (-1) - r \right] \frac{nmK}{(n+1)(m+1)} + r K
= 0 . \label{F(H(s,s)} 
\end{align}
It follows that $f(x,s) < 0$ for all $s > 0$ and $x \in
[0, H(s)]$, and (\ref{HJB-S}) has been established.

\underline{\em A probabilistic representation of $g$\/.}
Before addressing the proof of (\ref{HJB-W}), we first show
that, given any stopping time $\tau \leq \tau_\st$, where
$\tau_\st$ is defined by (\ref{opt-tau}),
\ben
g(x,s) = \e \left[ e^{-r \tau} g (X_\tau, S_\tau)
+ \int _0^\tau e^{-rt} f(X_t, S_t) \, dt + p \int _0^\tau
e^{-rt} S_t^{p-2} \, dS_t \right] . \label{g-FK}
\een
To this end, we assume that $(x,s) \in \W$ in what follows
without loss of generality.
Since the function $w(\cdot ,s)$ satisfies the ODE (\ref{ODE})
in the waiting region $\W$, we can see that
\ben
\half \sigma^2 x^2 g_{xx} (x,s) + \mu x g_x (x,s) -
rg(x,s) = - f(x,s) \quad \text{for all } s > 0 \text{ and }
x \in \mbox{} ]H(s), s[ . \label{Lg}
\een
Using It\^{o}'s formula, (\ref{HJB-BC}), the definition
of $g$ in (\ref{HJB-W}) and this calculation, we obtain
\begin{align}
g(x,s) = \mbox{} & e^{-r (\tau \wedge T)} g
(X_{\tau \wedge T}, S_{\tau \wedge T}) + \int
_0^{\tau \wedge T} e^{-rt} f(X_t, S_t) \, dt + p \int
_0^{\tau \wedge T} e^{-rt} S_t^{p-2} \, dS_t
\nonumber \\
& - \sigma \int _0^{\tau \wedge T} e^{-rt} g_x (X_t, S_t)
X_t \, dW_t . \nonumber
\end{align}
It follows that
\begin{align}
g(x,s) = \e \biggl[ & e^{-r (\tau \wedge T \wedge \tau_j)}
g (X_{\tau \wedge T \wedge \tau_j},
S_{\tau \wedge T \wedge \tau_j}) \nonumber \\
& + \int _0^{\tau \wedge T \wedge \tau_j} e^{-rt}
f(X_t, S_t) \, dt + p \int _0^{\tau \wedge T \wedge \tau_j}
e^{-rt} S_t^{p-2} \, dS_t \biggr] , \label{HJB-W-Ito}
\end{align}
where $(\tau_j)$ is a localising sequence of stopping
times for the stochastic integral.

Combining (\ref{w-growth}), (\ref{s/x-est}) and the
positivity of $w$ with the definition of $g$ in (\ref{HJB-W})
and the fact that $S$ is an increasing process, we can
see that
\be
\bigl| g (X_T, S_T) \bigr| \leq \left[ C + C S_T^\gamma
+ S_T^{p-1} + K \right] \quad \text{for all } T \leq \tau .
\ee
On the other hand, (\ref{s/x-est}), the definition of
$f$ in (\ref{HJB-S}), (\ref{f-deriv}) and the fact that 
$S$ is an increasing process imply that there exists
a constant $C_2 = C_2 (s) > 0$ such that
\be
\bigl| f (X_t, S_t) \bigr| \leq C_2 \bigl( 1 + S_t^{p-1}
\bigr) \quad \text{for all } t \leq \tau .
\ee
These estimates, the fact that $\gamma \in \mbox{}
]0,n[$, the assumption that $p-1 < n$ and
Lemma~\ref{supSgam-lem} imply that
\begin{align}
\e \left[ \sup _{T \geq 0} e^{-r (T \wedge \tau)}
\bigl| g (X_{T \wedge \tau}, S_{T \wedge \tau})
\bigr| \right] & \leq \e \left[ \sup _{T \geq 0}
e^{-r (T \wedge \tau)} \left[ C + C
S_{T \wedge \tau}^\gamma +
S_{T \wedge \tau}^{p-1} + K \right] \right]
\nonumber \\
& \leq \e \left[ \sup _{T \geq 0} e^{-rT} \left[
C + C S_T^\gamma + S_T^{p-1} + K \right] \right]
\nonumber \\
& < \infty \nonumber
\end{align}
and
\be
\e \left[ \int _0^\tau e^{-rt} \bigl| f(X_t, S_t) \bigr| \,
dt \right] \leq \e \left[ \int _0^\infty e^{-rt} C_2
\bigl( 1 + S_t^{p-1} \bigr) \, dt \right] < \infty .
\ee
In view of these observations and the fact that $S$
is an increasing process, we can pass to
the limits as $j \rightarrow \infty$ and
$T \rightarrow \infty$ in (\ref{HJB-W-Ito}) using
the dominated and the monotone convergence
theorems to obtain (\ref{g-FK}).

\underline{\em Proof of (\ref{HJB-W})\/.}
We first note that (\ref{a-quad}) and the definition
(\ref{Gamma12}) of $\Gamma_1$ imply that
\be
- f(s,s) = \left. \begin{cases} > 0 , & \text{if }
s^{p-1} > \Gamma_1^{-1} \\ < 0 , & \text{if }
s^{p-1} < \Gamma_1^{-1} \end{cases} \right\} .
\ee
Combining these inequalities with (\ref{f-deriv}) and
(\ref{F(H(s,s)}), we can see that
\be
- f(x,s) = \left. \begin{cases} > 0 , & \text{if }
s^{p-1} > \Gamma_1^{-1} \text{ and } x \in \mbox{}
]H(s), s[ \\ > 0 , & \text{if } s^{p-1} < \Gamma_1^{-1}
\text{ and } x \in \mbox{} ]H(s), \tilde{x} (s)[  \\ < 0 ,
& \text{if } s^{p-1} < \Gamma_1^{-1} \text{ and }
x \in \mbox{} ]\tilde{x} (s), s[ \end{cases} \right\} ,
\ee
where $\tilde{x} (s)$ is a unique point in $]H(s), s[$
for all $s>0$ such that $s^{p-1} < \Gamma_1^{-1}$.
In view of these inequalities, (\ref{Lg}) and the maximum
principle, we can see that, given any $s > 0$,
\begin{align}
\text{if } s^{p-1} \geq \Gamma_1^{-1} , & \text{ then the
function } g(\cdot ,s) \text{ has no positive maximum
inside } ]H(s), s[ , \label{MP1} \\
\text{if } s^{p-1} < \Gamma_1^{-1} , & \text{ then the
function } g(\cdot ,s) \text{ has no positive maximum
inside } ]H(s), \tilde{x} (s)[ , \label{MP2} \\
\text{if } s^{p-1} < \Gamma_1^{-1} , & \text{ then the
function } g(\cdot ,s) \text{ has no negative minimum
inside } ]\tilde{x} (s), s[ . \label{MP3}
\end{align}

To proceed further, we use the identity
\begin{align}
g_{xx} (x,s) = \mbox{} & n(n-1) \frac{-(m+1) s^p H^{-1}(s) + mK}
{n-m} H^{-n}(s) x^{n-2} \nonumber \\
& + m(m-1) \frac{(n+1) s^p H^{-1}(s) - nK}{n-m} H^{-m}(s)
x^{m-2} - 2 s^p x^{-3} , \nonumber
\end{align}
which holds true in $\W$ by the definition (\ref{w}) of $w$,
to calculate
\begin{align}
\lim_{x \downarrow H(s)} g_{xx} (x,s) & = - \bigl[ 1 + n+m
+ nm \bigr] s^p H^{-1}(s) + nmK H^{-2}(s) \nonumber \\
& \stackrel{(\ref{a-quad})}{=} - \frac{2}{\sigma^2} f \bigl(
H(s), s \bigr)
H^{-2} (s) \stackrel{(\ref{F(H(s,s)})}{>} 0 . \nonumber
\end{align}
This result and the identities $g \bigl(H(s), s \bigr) = g_x
\bigl( H(s), s \bigr) = 0$, which follow from the $C^1$-continuity
of $w(\cdot ,s)$ at $H(s)$, imply that
\ben
g_x \bigl( H(s) + \varepsilon , s \bigr) > 0 \quad \text{and}
\quad g \bigl( H(s) + \varepsilon , s \bigr) > 0 \quad
\text{for all } \varepsilon > 0 \text{ sufficiently small} .
\label{g-at_H}
\een
Combining this observation with (\ref{MP1}) we obtain
(\ref{HJB-W}) for all $s>0$ such that $s^{p-1} \geq
\Gamma_1^{-1}$ and $x \in \mbox{} ]H(s), s[$.
On the other hand, combining (\ref{g-at_H}) with
(\ref{MP2})--(\ref{MP3}), we obtain (\ref{HJB-W})
for all $s>0$ such that $s^{p-1} \leq K$ and
$x \in \mbox{} ]H(s), s[$ because $g(s,s) \geq 0$ if
$s^{p-1} \leq K$ thanks to the positivity of $w$.
It follows that
\ben
g(x,s) = w(x,s) - \frac{s^p}{x} + K \geq 0 \quad \text{if }
s^{p-1} \in \mbox{} ]0, K] \cup [\Gamma_1^{-1}, \infty[
\text{ and } x \in \mbox{} ]H(s), s[ . \label{HJB-W-C0}
\een
In particular, (\ref{HJB-W}) holds true if
\ben
\Gamma_1^{-1} \leq K  \quad \Leftrightarrow \quad
\Gamma_1 \geq \Gamma_2 \quad \Leftrightarrow \quad
\mu \geq \sigma^2 \label{G1<K}
\een
(see also (\ref{Gamma}) and (\ref{Gamma12})).

To establish (\ref{HJB-W}) if the problem data is such
that (\ref{G1<K}) is not true, we argue by contradiction.
In view of (\ref{MP2})--(\ref{MP3}) and (\ref{HJB-W-C0}),
we therefore assume that $K < \Gamma_1^{-1}$ and that
there exist strictly positive $\tilde{s}_1 < \tilde{s}_2$
such that $\tilde{s}_1^{p-1} , \tilde{s}_2^{p-1} \in
[K, \Gamma_1^{-1}]$,
\ben
g(x,s) \left. \begin{cases} < 0 & \text{for all } x = s \in
\mbox{} ]\tilde{s}_1, \tilde{s}_2[ \\ \geq 0 & \text{for all }
s \in \mbox{} ]0, \tilde{s}_1] \cup [\tilde{s}_2, \infty[
\text{ and } x \in \mbox{} ]H(s), s[ \end{cases} \right\}
. \label{HJB-W-C1}
\een
Also, we note that (\ref{g-at_H}) implies that there
exists $\tilde{\varepsilon} > 0$ such that
\ben
H(\tilde{s}_2) < \tilde{s}_2 - \tilde{\varepsilon}
\quad \text{and} \quad
g \bigl( H(\tilde{s}_2) , s \bigr) > 0 \quad
\text{for all } s \in [\tilde{s}_2 - \tilde{\varepsilon},
\tilde{s}_2[ . \label{HJB-W-C2}
\een
Given such an $\tilde{\varepsilon} > 0$ fixed, we
consider the solution to (\ref{X}) with initial condition
$X_0 = \tilde{s}_2 - \tilde{\varepsilon}$ and the
running maximum process $S$ given by (\ref{S})
with initial condition $S_0 = \tilde{s}_2 -
\tilde{\varepsilon}$.
Also, we define
\be
\bar{S}_t = \tilde{s}_2 \vee S_t , \quad \text{for }
t \geq 0 , \quad \text{and} \quad
\tau = \inf \bigl\{ t \geq 0 \mid \ X_t = H(\tilde{s}_2)
\vee H(S_t) \bigr\} .
\ee
Using (\ref{g-FK}), (\ref{HJB-W-C1})--(\ref{HJB-W-C2}),
the identity $g \bigl(H(s), s \bigr) = 0$ that holds true
for all $s>0$ and the fact that the function $f(x,\cdot)
: [x,\infty[ \mbox{} \rightarrow \R$ is strictly decreasing
for all $x>0$, which follows from the calculation 
\be
f_s (x,s) = p \left[ \half \sigma^2 (-1)^2 + \left( \mu - \half
\sigma^2 \right) (-1) - r \right] \frac{s^{p-1}}{x}
\stackrel{(\ref{a-quad})}{<} 0 \quad
\text{for all } s > 0 \text{ and } x \in \mbox{} ]0, s[ ,
\ee
we obtain
\begin{align}
0 > \mbox{} & g(\tilde{s}_2 - \tilde{\varepsilon},
\tilde{s}_2 - \tilde{\varepsilon}) \nonumber \\
= \mbox{} & \e \left[ e^{-r \tau} g (X_\tau, S_\tau)
+ \int _0^\tau e^{-rt} f(X_t, S_t) \, dt + p \int _0^\tau
e^{-rt} S_t^{p-2} \, dS_t \right] \nonumber \\
= \mbox{} & \e \biggl[ e^{-r \tau} g \bigl( H(\tilde{s}_2),
S_\tau \bigr) {\bf 1} _{\{ S_\tau < \tilde{s}_2 \}}
+ \int _0^\tau e^{-rt} f(X_t, S_t) \, dt \nonumber \\
& \hspace{5mm} \mbox{} + p \int _0^\tau e^{-rt}
{\bf 1} _{\{ S_t < \tilde{s}_2 \}} S_t^{p-2} \, dS_t
+ p \int _0^\tau e^{-rt} \bar{S}_t^{p-2} \, d\bar{S}_t
\biggr] \nonumber \\
> \mbox{} & \e \left[ \int _0^\tau e^{-rt} f(X_t, \bar{S}_t)
\, dt + p \int _0^\tau e^{-rt} \bar{S}_t^{p-2}
\, d\bar{S}_t \right] \nonumber \\
= \mbox{} & \e \left[ e^{-r \tau} g \bigl( X_\tau,
\bar{S}_\tau \bigr) + \int _0^\tau e^{-rt}
f(X_t, \bar{S}_t) \, dt + p \int _0^\tau e^{-rt}
\bar{S}_t^{p-2} \, d\bar{S}_t \right] \nonumber \\
= \mbox{} & g(\tilde{s}_2 - \tilde{\varepsilon},
\tilde{s}_2) \nonumber \\
\geq \mbox{} & 0 , \nonumber
\end{align}
which is a contradiction. 
\mbox{}\hfill$\Box$


\newpage
\begin{picture}(160,110)
\put(20,15){\begin{picture}(120,95) 

\put(0,0){\line(1,0){120}}
\put(120,0){\line(0,1){95}}
\put(0,0){\line(0,1){95}}
\put(0,95){\line(1,0){120}}

\put(10,10){\vector(1,0){100}}
\put(10,10){\vector(0,1){80}}
\put(6,88){$x$}
\put(108,6){$s$}

\put(0,0){\qbezier(10,10)(60,50)(110,90)}

\put(0,-0.2){\qbezier(10,10)(25,11.7)(50,30)}
\put(0,-0.1){\qbezier(10,10)(25,11.7)(50,30)}
\put(0,0){\qbezier(10,10)(25,11.7)(50,30)}
\put(0,0.1){\qbezier(10,10)(25,11.7)(50,30)}
\put(0,0.2){\qbezier(10,10)(25,11.7)(50,30)}

\put(0,-0.2){\qbezier(50,30)(85,55)(110,61.5)}
\put(0,-0.1){\qbezier(50,30)(85,55)(110,61.5)}
\put(0,0){\qbezier(50,30)(85,55)(110,61.5)}
\put(0,0.1){\qbezier(50,30)(85,55)(110,61.5)}
\put(0,0.2){\qbezier(50,30)(85,55)(110,61.5)}

\put(102,70){$H(s)$}
\put(105,68){\vector(0,-1){7}}

\put(70,50){$\mathcal W$}
\put(85,28){$\mathcal S$}

\end{picture}}
\end{picture}

\mbox{} \vspace{-15mm}

\noindent
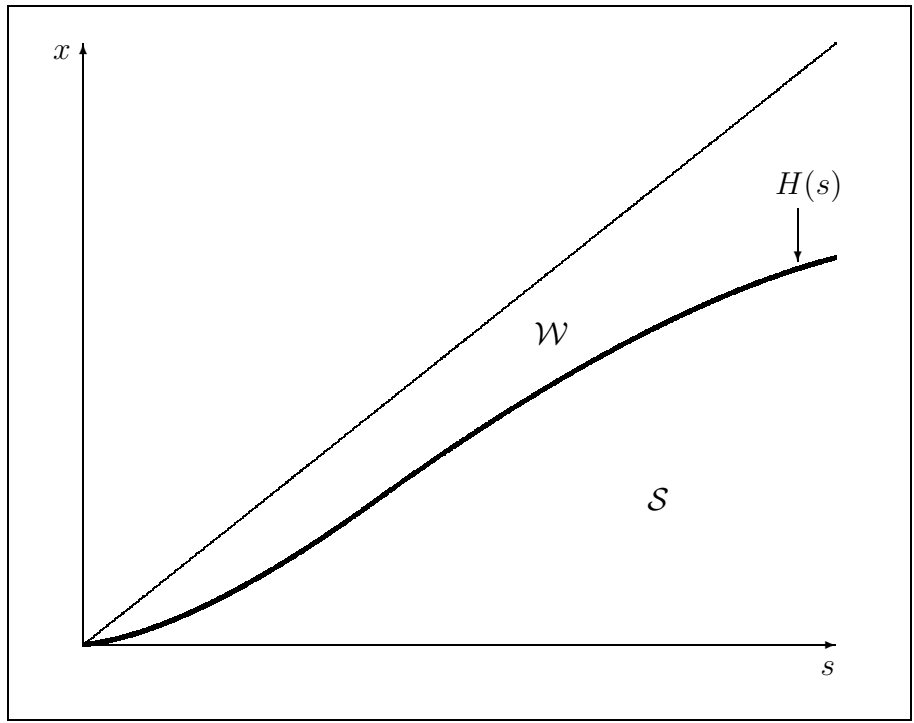
{\bf Figure 1} Depiction of the
free-boundary function $H$ separating the stopping
region $\st$ from the waiting region $\W$.

\newpage
\begin{picture}(160,110)
\put(20,15){\begin{picture}(120,95) 

\put(0,0){\line(1,0){120}}
\put(120,0){\line(0,1){95}}
\put(0,0){\line(0,1){95}}
\put(0,95){\line(1,0){120}}

\put(10,10){\vector(1,0){100}}
\put(10,10){\vector(0,1){80}}
\put(6,88){$x$}
\put(108,6){$s$}

\put(0,0){\qbezier(10,10)(60,50)(110,90)}

\put(0,0){\qbezier(10,10)(30,50)(110,75)}

\color{olive}
\put(0,-0.1){\qbezier(10,10)(60,50)(75,62)}
\put(0,0){\qbezier(10,10)(60,50)(75,62)}
\put(0,0.1){\qbezier(10,10)(60,50)(75,62)}

\put(0,-0.1){\qbezier(75,62)(90,68.8)(110,75)}
\put(0,0){\qbezier(75,62)(90,68.8)(110,75)}
\put(0,0.1){\qbezier(75,62)(90,68.8)(110,75)}

\put(10,9.9){\line(1,0){98.5}}
\put(10,10){\line(1,0){98.5}}
\put(10,10.1){\line(1,0){98.5}}
\color{black}

\put(0,0){\qbezier(10,10)(31,45)(110,62)}

\put(9,30){\line(1,0){2}}
\put(0,0){\qbezier[100](10,30)(60,30)(110,30)}
\put(6,29.2){$\delta$}

\color{red}
\put(0,-0.2){\qbezier(10,10)(25,11.7)(50,30)}
\put(0,-0.1){\qbezier(10,10)(25,11.7)(50,30)}
\put(0,0){\qbezier(10,10)(25,11.7)(50,30)}
\put(0,0.1){\qbezier(10,10)(25,11.7)(50,30)}
\put(0,0.2){\qbezier(10,10)(25,11.7)(50,30)}

\put(0,-0.2){\qbezier(50,30)(85,55)(110,61.5)}
\put(0,-0.1){\qbezier(50,30)(85,55)(110,61.5)}
\put(0,0){\qbezier(50,30)(85,55)(110,61.5)}
\put(0,0.1){\qbezier(50,30)(85,55)(110,61.5)}
\put(0,0.2){\qbezier(50,30)(85,55)(110,61.5)}

\color{cyan}
\put(0,-0.2){\qbezier(10,10)(25,17)(43,30)}
\put(0,-0.1){\qbezier(10,10)(25,17)(43,30)}
\put(0,0){\qbezier(10,10)(25,17)(43,30)}
\put(0,0.1){\qbezier(10,10)(25,17)(43,30)}
\put(0,0.2){\qbezier(10,10)(25,17)(43,30)}

\put(0,-0.2){\qbezier(43,30)(85,60)(95,69.2)}
\put(0,-0.1){\qbezier(43,30)(85,60)(95,69.2)}
\put(0,0){\qbezier(43,30)(85,60)(95,69.2)}
\put(0,0.1){\qbezier(43,30)(85,60)(95,69.2)}
\put(0,0.2){\qbezier(43,30)(85,60)(95,69.2)}

\put(0,-0.2){\qbezier(10,10)(35,11)(60,30)}
\put(0,-0.1){\qbezier(10,10)(35,11)(60,30)}
\put(0,0){\qbezier(10,10)(35,11)(60,30)}
\put(0,0.1){\qbezier(10,10)(35,11)(60,30)}
\put(0,0.2){\qbezier(10,10)(35,11)(60,30)}

\put(0,-0.2){\qbezier(60,30)(85,47)(110,48)}
\put(0,-0.1){\qbezier(60,30)(85,47)(110,48)}
\put(0,0){\qbezier(60,30)(85,47)(110,48)}
\put(0,0.1){\qbezier(60,30)(85,47)(110,48)}
\put(0,0.2){\qbezier(60,30)(85,47)(110,48)}
\color{black}

\put(35,9){\line(0,1){2}}
\put(43,9){\line(0,1){2}}
\put(50,9){\line(0,1){2}}
\put(60,9){\line(0,1){2}}

\put(0,0){\qbezier[20](35,10)(35,20)(35,30)}
\put(0,0){\qbezier[20](43,10)(43,20)(43,30)}
\put(0,0){\qbezier[20](50,10)(50,20)(50,30)}
\put(0,0){\qbezier[20](60,10)(60,20)(60,30)}

\put(34,6){$s_\dagger$}
\put(42,6){$s_*^1$}
\put(49,6){$s_\circ$}
\put(59,6){$s_*^2$}

\put(96,65){$H(s;s_*^1)$}
\put(95.7,66.3){\vector(-1,0){3.4}}

\put(96,52){$H(s;s_\circ)$}
\put(95.7,53.3){\vector(-1,0){6}}

\put(96,40){$H(s;s_*^2)$}
\put(95.7,41.3){\vector(-1,0){14}}

\put(23,43){$c s^p$}
\put(25.7,41.5){\vector(1,-1){7.7}}

\put(63,76){$s \wedge [\Gamma s^p]$}
\put(64.3,74){\vector(0,-1){20}}
\put(75,74){\vector(1,-1){8}}

\end{picture}}
\end{picture}

\mbox{} \vspace{-15mm}

\noindent
{\bf Figure 2 ($\pmb{p<1}$)}
Illustration of Lemma~\ref{H-lem1}.(I) for $s_\circ (\delta)
= s^\circ (\delta)$.
The free-boundary function $H(\cdot) = H(\cdot; s_\circ)$
that separates the stopping region $\st$ from the waiting
region $\W$ is plotted in red.
The intersection of $\R^2$ with the boundary of the
domain $\domH$ in which we consider solutions to
the ODE (\ref{H-ODE}) satisfying (\ref{H-IC}) is designated
by green.
Every solution to the ODE (\ref{H-ODE}) satisfying
(\ref{H-IC}) with $s_* \in \mbox{} ]s_\dagger, s_\circ[$
(resp., $s_* > s_\circ$) hits the upper part of the boundary
of $\domH$ in the picture (resp., has asymptotic growth
as $s \rightarrow \infty$ that is of different order than the
one of the free-boundary): such solutions are plotted
in blue.

\newpage
\begin{picture}(160,110)
\put(20,15){\begin{picture}(120,95) 

\put(0,0){\line(1,0){120}}
\put(120,0){\line(0,1){95}}
\put(0,0){\line(0,1){95}}
\put(0,95){\line(1,0){120}}

\put(10,10){\vector(1,0){100}}
\put(10,10){\vector(0,1){80}}
\put(6,88){$x$}
\put(108,6){$s$}

\put(0,0){\qbezier(10,10)(60,50)(110,90)}

\put(0,0){\qbezier(10,10)(60,20)(80,90)}

\color{olive}
\put(0,-0.1){\qbezier(10,10)(45,16.7)(66.1,54.8)}
\put(0,0){\qbezier(10,10)(45,16.7)(66.1,54.8)}
\put(0,0.1){\qbezier(10,10)(45,16.7)(66.1,54.8)}

\put(0,-0.1){\qbezier(66.1,54.8)(88.05,72.4)(110,90)}
\put(0,0){\qbezier(66.1,54.8)(88.05,72.4)(110,90)}
\put(0,0.1){\qbezier(66.1,54.8)(88.05,72.4)(110,90)}

\put(10,9.9){\line(1,0){98.5}}
\put(10,10){\line(1,0){98.5}}
\put(10,10.1){\line(1,0){98.5}}
\color{black}

\put(0,0){\qbezier(10,10)(60,35)(110,60)}

\put(9,30){\line(1,0){2}}
\put(0,0){\qbezier[100](10,30)(60,30)(110,30)}
\put(6,29.2){$\delta$}

\color{red}
\put(0,-0.2){\qbezier(10,10)(46,11)(70,30)}
\put(0,-0.1){\qbezier(10,10)(46,11)(70,30)}
\put(0,0){\qbezier(10,10)(46,11)(70,30)}
\put(0,0.1){\qbezier(10,10)(46,11)(70,30)}
\put(0,0.2){\qbezier(10,10)(46,11)(70,30)}

\put(0,-0.2){\qbezier(70,30)(90,47)(110,59.6)}
\put(0,-0.1){\qbezier(70,30)(90,47)(110,59.6)}
\put(0,0){\qbezier(70,30)(90,47)(110,59.6)}
\put(0,0.1){\qbezier(70,30)(90,47)(110,59.6)}
\put(0,0.2){\qbezier(70,30)(90,47)(110,59.6)}

\color{cyan}
\put(0,-0.2){\qbezier(10,10)(40,11)(60,30)}
\put(0,-0.1){\qbezier(10,10)(40,11)(60,30)}
\put(0,0){\qbezier(10,10)(40,11)(60,30)}
\put(0,0.1){\qbezier(10,10)(40,11)(60,30)}
\put(0,0.2){\qbezier(10,10)(40,11)(60,30)}

\put(0,-0.2){\qbezier(60,30)(80,50)(86,70)}
\put(0,-0.1){\qbezier(60,30)(80,50)(86,70)}
\put(0,0){\qbezier(60,30)(80,50)(86,70)}
\put(0,0.1){\qbezier(60,30)(80,50)(86,70)}
\put(0,0.2){\qbezier(60,30)(80,50)(86,70)}

\put(0,-0.2){\qbezier(10,10)(55,11)(80,30)}
\put(0,-0.1){\qbezier(10,10)(55,11)(80,30)}
\put(0,0){\qbezier(10,10)(55,11)(80,30)}
\put(0,0.1){\qbezier(10,10)(55,11)(80,30)}
\put(0,0.2){\qbezier(10,10)(55,11)(80,30)}

\put(0,-0.2){\qbezier(80,30)(100,45)(110,48)}
\put(0,-0.1){\qbezier(80,30)(100,45)(110,48)}
\put(0,0){\qbezier(80,30)(100,45)(110,48)}
\put(0,0.1){\qbezier(80,30)(100,45)(110,48)}
\put(0,0.2){\qbezier(80,30)(100,45)(110,48)}
\color{black}

\put(47.7,9){\line(0,1){2}}
\put(60,9){\line(0,1){2}}
\put(70,9){\line(0,1){2}}
\put(80,9){\line(0,1){2}}

\put(0,0){\qbezier[20](47.7,10)(47.7,20)(47.7,30)}
\put(0,0){\qbezier[20](60,10)(60,20)(60,30)}
\put(0,0){\qbezier[20](70,10)(70,20)(70,30)}
\put(0,0){\qbezier[20](80,10)(80,20)(80,30)}

\put(46.7,6){$s_\dagger$}
\put(59,6){$s_*^1$}
\put(69,6){$s_\circ$}
\put(79,6){$s_*^2$}

\put(96,70){$H(s;s_\circ)$}
\put(105,68){\vector(0,-1){11}}

\put(86,60){$H(s;s_*^1)$}
\put(85.55,61.3){\vector(-1,0){2.3}}

\put(96,35){$H(s;s_*^2)$}
\put(105,39){\vector(0,1){6.5}}

\put(16,43){$c s$}
\put(18.7,41.5){\vector(1,-1){17.7}}

\put(42,63){$[\Gamma s^p] \wedge s$}
\put(51,60.2){\vector(1,-1){11.2}}
\put(61.3,64){\vector(1,0){15.5}}

\end{picture}}
\end{picture}

\mbox{} \vspace{-15mm}

\noindent
{\bf Figure 3 ($\pmb{p>1}$)}
Illustration of Lemma~\ref{H-lem1}.(II) for $s_\circ (\delta)
= s^\circ (\delta)$.
The free-boundary function $H(\cdot) = H(\cdot; s_\circ)$
that separates the stopping region $\st$ from the waiting
region $\W$ is plotted in red.
The intersection of $\R^2$ with the boundary of the
domain $\domH$ in which we consider solutions to
the ODE (\ref{H-ODE}) satisfying (\ref{H-IC}) is designated
by green.
Every solution to the ODE (\ref{H-ODE}) satisfying
(\ref{H-IC}) with $s_* \in \mbox{} ]s_\dagger, s_\circ[$
(resp., $s_* > s_\circ$) hits the upper part of the boundary
of $\domH$ in the picture (resp., has asymptotic growth
as $s \rightarrow \infty$ that is of different order than the
one of the free-boundary): such solutions are plotted
in blue.

\newpage
\begin{picture}(160,110)
\put(20,15){\begin{picture}(120,95) 

\put(0,0){\line(1,0){120}}
\put(120,0){\line(0,1){95}}
\put(0,0){\line(0,1){95}}
\put(0,95){\line(1,0){120}}

\put(10,10){\vector(1,0){100}}
\put(10,10){\vector(0,1){80}}
\put(6,88){$x$}
\put(108,6){$s$}

\put(9,65){\line(1,0){2}}
\put(6,64.2){$\Gamma$}

\put(9,50){\line(1,0){2}}
\put(0,0){\qbezier(10,50)(60,50)(110,50)}
\put(6,49.2){$c$}

\put(0,0){\qbezier(10,10)(20,65)(80,90)}
\put(0,0){\qbezier(10,65)(60,65)(110,65)}
\put(70.5,75){$s^{1-p}$}
\put(69.5,76){\vector(-1,0){15}}

\put(0,0){\qbezier[150](11,88)(15,15)(108,11)}
\put(21,75){$\delta s^{-p}$}
\put(19.5,76){\vector(-1,0){6.5}}

\color{olive}
\put(0,-0.1){\qbezier(10,10)(16,43)(39.8,65)}
\put(0,0){\qbezier(10,10)(16,43)(39.8,65)}
\put(0,0.1){\qbezier(10,10)(16,43)(39.8,65)}

\put(0,-0.1){\qbezier(39.8,65)(80,65)(110,65)}
\put(0,0){\qbezier(39.8,65)(80,65)(110,65)}
\put(0,0.1){\qbezier(39.8,65)(80,65)(110,65)}

\put(10,9.9){\line(1,0){98.5}}
\put(10,10){\line(1,0){98.5}}
\put(10,10.1){\line(1,0){98.5}}
\color{black}

\color{red}
\put(0,-0.2){\qbezier(10,10)(25,21)(40,30)}
\put(0,-0.1){\qbezier(10,10)(25,21)(40,30)}
\put(0,0){\qbezier(10,10)(25,21)(40,30)}
\put(0,0.1){\qbezier(10,10)(25,21)(40,30)}
\put(0,0.2){\qbezier(10,10)(25,21)(40,30)}

\put(0,-0.2){\qbezier(40,30)(70,48)(110,49.5)}
\put(0,-0.1){\qbezier(40,30)(70,48)(110,49.5)}
\put(0,0){\qbezier(40,30)(70,48)(110,49.5)}
\put(0,0.1){\qbezier(40,30)(70,48)(110,49.5)}
\put(0,0.2){\qbezier(40,30)(70,48)(110,49.5)}

\color{blue}
\put(0,-0.2){\qbezier(10,10)(25,20)(42,29)}
\put(0,-0.1){\qbezier(10,10)(25,20)(42,29)}
\put(0,0){\qbezier(10,10)(25,20)(42,29)}
\put(0,0.1){\qbezier(10,10)(25,20)(42,29)}
\put(0,0.2){\qbezier(10,10)(25,20)(42,29)}

\put(0,-0.1){\qbezier(42,29)(70,44)(110,49)}
\put(0,-0.1){\qbezier(42,29)(70,44)(110,49)}
\put(0,0){\qbezier(42,29)(70,44)(110,49)}
\put(0,0.1){\qbezier(42,29)(70,44)(110,49)}
\put(0,0.2){\qbezier(42,29)(70,44)(110,49)}

\color{cyan}
\put(0,-0.2){\qbezier(10,10)(27,26)(35,34)}
\put(0,-0.1){\qbezier(10,10)(27,26)(35,34)}
\put(0,0){\qbezier(10,10)(27,26)(35,34)}
\put(0,0.1){\qbezier(10,10)(27,26)(35,34)}
\put(0,0.2){\qbezier(10,10)(27,26)(35,34)}

\put(0,-0.2){\qbezier(35,34)(57,55)(80,64.5)}
\put(0,-0.1){\qbezier(35,34)(57,55)(80,64.5)}
\put(0,0){\qbezier(35,34)(57,55)(80,64.5)}
\put(0,0.1){\qbezier(35,34)(57,55)(80,64.5)}
\put(0,0.2){\qbezier(35,34)(57,55)(80,64.5)}

\put(0,-0.2){\qbezier(10,10)(37,26.4)(45,27)}
\put(0,-0.1){\qbezier(10,10)(37,26.4)(45,27)}
\put(0,0){\qbezier(10,10)(37,26.4)(45,27)}
\put(0,0.1){\qbezier(10,10)(37,26.4)(45,27)}
\put(0,0.2){\qbezier(10,10)(37,26.4)(45,27)}

\put(0,-0.2){\qbezier(45,27)(48,26.9)(55,24)}
\put(0,-0.1){\qbezier(45,27)(48,26.9)(55,24)}
\put(0,0){\qbezier(45,27)(48,26.9)(55,24)}
\put(0,0.1){\qbezier(45,27)(48,26.9)(55,24)}
\put(0,0.2){\qbezier(45,27)(48,26.9)(55,24)}

\put(0,-0.2){\qbezier(55,24)(75,15)(108,11.6)}
\put(0,-0.1){\qbezier(55,24)(75,15)(108,11.6)}
\put(0,0){\qbezier(55,24)(75,15)(108,11.6)}
\put(0,0.1){\qbezier(55,24)(75,15)(108,11.6)}
\put(0,0.2){\qbezier(55,24)(75,15)(108,11.6)}
\color{black}

\put(24,5){\line(0,1){6}}
\put(35,9){\line(0,1){2}}
\put(40,9){\line(0,1){2}}
\put(42,4){\line(0,1){7}}
\put(45,9){\line(0,1){2}}

\put(0,0){\qbezier[40](24,10)(24,28)(24,46)}
\put(0,0){\qbezier[30](35,10)(35,22)(35,34)}
\put(0,0){\qbezier[20](40,10)(40,20)(40,30)}
\put(0,0){\qbezier[20](42,10)(42,20)(42,29)}
\put(0,0){\qbezier[20](45,10)(45,20)(45,27)}

\put(23,2){$s_\dagger$}
\put(33.5,6){$s_*^1$}
\put(38,6){$s_\circ$}
\put(41,1){$s^\dagger$}
\put(44,6){$s_*^2$}

\put(96,38.5){$h(s;s_\circ)$}
\put(102,42.5){\vector(0,1){6.2}}

\put(78,56){$h(s;s_*^1)$}
\put(77,57.3){\vector(-1,0){10}}

\put(96,20){$h(s;s_*^2)$}
\put(102,18){\vector(0,-1){5}}

\put(72,30){$\mathfrak{s}^{\{ -1 \}} (s)$}
\put(79,35){\vector(0,1){7.5}}

\end{picture}}
\end{picture}

\mbox{} \vspace{-15mm}

\noindent
{\bf Figure 4 ($\pmb{p<1}$)}
Illustration of the proof of Lemma~\ref{H-lem1}.(I) for
$s_\circ (\delta) = s^\circ (\delta)$.
The identity $h(s;s_*) = s^{-p} H(s;s_*)$ for all $s > 0$
relates the solutions to (\ref{h-ODE1}) for
$s_* = s_*^1, s_\circ, s_*^2$ plotted here with the
solutions to the ODE (\ref{H-ODE}) satisfying (\ref{H-IC})
that are plotted in Figure~2.
Furthermore, the intersection of $\R^2$ with the boundary
of the domain $\domh^1$ in which we consider solutions to
(\ref{h-ODE1}) is designated by green.

\newpage
\begin{picture}(160,110)
\put(20,15){\begin{picture}(120,95) 

\put(0,0){\line(1,0){120}}
\put(120,0){\line(0,1){95}}
\put(0,0){\line(0,1){95}}
\put(0,95){\line(1,0){120}}

\put(10,10){\vector(1,0){100}}
\put(10,10){\vector(0,1){80}}
\put(6,88){$x$}
\put(108,6){$s$}

\put(9,65){\line(1,0){2}}
\put(6,64.2){1}

\put(9,50){\line(1,0){2}}
\put(0,0){\qbezier(10,50)(60,50)(110,50)}
\put(6,49.2){$c$}

\put(0,0){\qbezier(10,10)(20,65)(80,90)}
\put(0,0){\qbezier(10,65)(60,65)(110,65)}
\put(70,75){$\Gamma s^{p-1}$}
\put(69.5,76){\vector(-1,0){15}}

\put(0,0){\qbezier[150](11,88)(15,15)(108,11)}
\put(21,75){$\delta s^{-1}$}
\put(19.5,76){\vector(-1,0){6.5}}

\color{olive}
\put(0,-0.1){\qbezier(10,10)(16,43)(39.8,65)}
\put(0,0){\qbezier(10,10)(16,43)(39.8,65)}
\put(0,0.1){\qbezier(10,10)(16,43)(39.8,65)}

\put(0,-0.1){\qbezier(39.8,65)(80,65)(110,65)}
\put(0,0){\qbezier(39.8,65)(80,65)(110,65)}
\put(0,0.1){\qbezier(39.8,65)(80,65)(110,65)}

\put(10,9.9){\line(1,0){98.5}}
\put(10,10){\line(1,0){98.5}}
\put(10,10.1){\line(1,0){98.5}}
\color{black}

\color{red}
\put(0,-0.2){\qbezier(10,10)(25,21)(40,30)}
\put(0,-0.1){\qbezier(10,10)(25,21)(40,30)}
\put(0,0){\qbezier(10,10)(25,21)(40,30)}
\put(0,0.1){\qbezier(10,10)(25,21)(40,30)}
\put(0,0.2){\qbezier(10,10)(25,21)(40,30)}

\put(0,-0.2){\qbezier(40,30)(70,48)(110,49.5)}
\put(0,-0.1){\qbezier(40,30)(70,48)(110,49.5)}
\put(0,0){\qbezier(40,30)(70,48)(110,49.5)}
\put(0,0.1){\qbezier(40,30)(70,48)(110,49.5)}
\put(0,0.2){\qbezier(40,30)(70,48)(110,49.5)}

\color{blue}
\put(0,-0.2){\qbezier(10,10)(25,20)(42,29)}
\put(0,-0.1){\qbezier(10,10)(25,20)(42,29)}
\put(0,0){\qbezier(10,10)(25,20)(42,29)}
\put(0,0.1){\qbezier(10,10)(25,20)(42,29)}
\put(0,0.2){\qbezier(10,10)(25,20)(42,29)}

\put(0,-0.1){\qbezier(42,29)(70,44)(110,49)}
\put(0,-0.1){\qbezier(42,29)(70,44)(110,49)}
\put(0,0){\qbezier(42,29)(70,44)(110,49)}
\put(0,0.1){\qbezier(42,29)(70,44)(110,49)}
\put(0,0.2){\qbezier(42,29)(70,44)(110,49)}

\color{cyan}
\put(0,-0.2){\qbezier(10,10)(27,26)(35,34)}
\put(0,-0.1){\qbezier(10,10)(27,26)(35,34)}
\put(0,0){\qbezier(10,10)(27,26)(35,34)}
\put(0,0.1){\qbezier(10,10)(27,26)(35,34)}
\put(0,0.2){\qbezier(10,10)(27,26)(35,34)}

\put(0,-0.2){\qbezier(35,34)(57,55)(80,64.5)}
\put(0,-0.1){\qbezier(35,34)(57,55)(80,64.5)}
\put(0,0){\qbezier(35,34)(57,55)(80,64.5)}
\put(0,0.1){\qbezier(35,34)(57,55)(80,64.5)}
\put(0,0.2){\qbezier(35,34)(57,55)(80,64.5)}

\put(0,-0.2){\qbezier(10,10)(37,26.4)(45,27)}
\put(0,-0.1){\qbezier(10,10)(37,26.4)(45,27)}
\put(0,0){\qbezier(10,10)(37,26.4)(45,27)}
\put(0,0.1){\qbezier(10,10)(37,26.4)(45,27)}
\put(0,0.2){\qbezier(10,10)(37,26.4)(45,27)}

\put(0,-0.2){\qbezier(45,27)(48,26.9)(55,24)}
\put(0,-0.1){\qbezier(45,27)(48,26.9)(55,24)}
\put(0,0){\qbezier(45,27)(48,26.9)(55,24)}
\put(0,0.1){\qbezier(45,27)(48,26.9)(55,24)}
\put(0,0.2){\qbezier(45,27)(48,26.9)(55,24)}

\put(0,-0.2){\qbezier(55,24)(75,15)(108,11.6)}
\put(0,-0.1){\qbezier(55,24)(75,15)(108,11.6)}
\put(0,0){\qbezier(55,24)(75,15)(108,11.6)}
\put(0,0.1){\qbezier(55,24)(75,15)(108,11.6)}
\put(0,0.2){\qbezier(55,24)(75,15)(108,11.6)}
\color{black}

\put(24,5){\line(0,1){6}}
\put(35,9){\line(0,1){2}}
\put(40,9){\line(0,1){2}}
\put(42,4){\line(0,1){7}}
\put(45,9){\line(0,1){2}}

\put(0,0){\qbezier[40](24,10)(24,28)(24,46)}
\put(0,0){\qbezier[30](35,10)(35,22)(35,34)}
\put(0,0){\qbezier[20](40,10)(40,20)(40,30)}
\put(0,0){\qbezier[20](42,10)(42,20)(42,29)}
\put(0,0){\qbezier[20](45,10)(45,20)(45,27)}

\put(23,2){$s_\dagger$}
\put(33.5,6){$s_*^1$}
\put(38,6){$s_\circ$}
\put(41,1){$s^\dagger$}
\put(44,6){$s_*^2$}

\put(96,38.5){$h(s;s_\circ)$}
\put(102,42.5){\vector(0,1){6.2}}

\put(78,56){$h(s;s_*^1)$}
\put(77,57.3){\vector(-1,0){10}}

\put(96,20){$h(s;s_*^2)$}
\put(102,18){\vector(0,-1){5}}

\put(72,30){$\mathfrak{s}^{\{ -1 \}} (s)$}
\put(79,35){\vector(0,1){7.5}}

\end{picture}}
\end{picture}

\mbox{} \vspace{-15mm}

\noindent
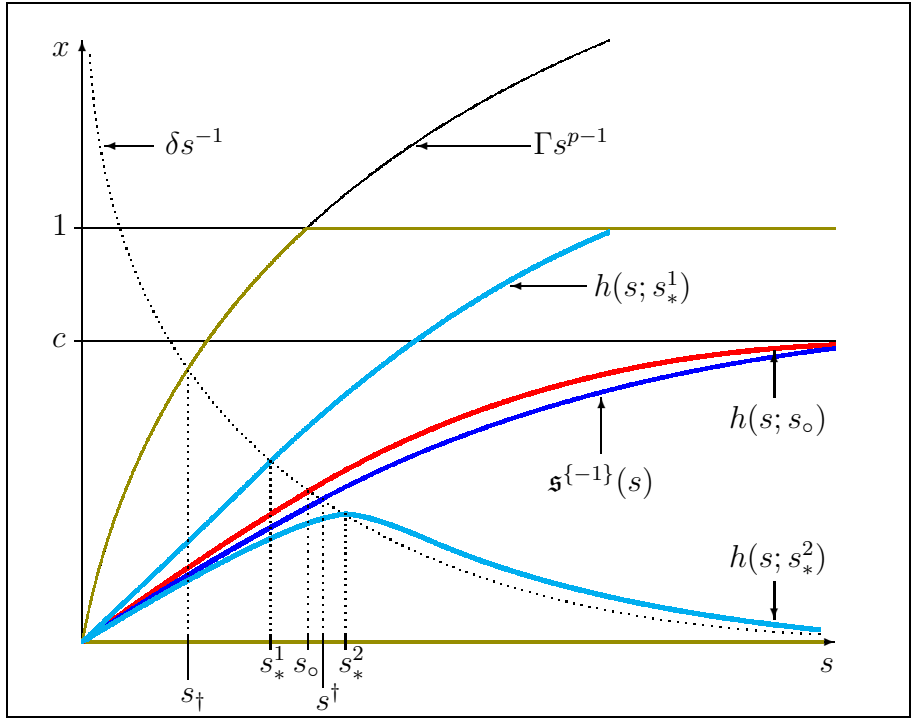
{\bf Figure 5 ($\pmb{p>1}$)}
Illustration of the proof of Lemma~\ref{H-lem1}.(II) for
$s_\circ (\delta) = s^\circ (\delta)$.
The identity $h(s;s_*) = s^{-1} H(s;s_*)$ for all $s > 0$
relates the solutions to (\ref{h-ODE2}) for
$s_* = s_*^1, s_\circ, s_*^2$ plotted here with the
solutions to the ODE (\ref{H-ODE}) satisfying (\ref{H-IC})
that are plotted in Figure~3.
Furthermore, the intersection of $\R^2$ with the boundary
of the domain $\domh^2$ in which we consider solutions to
(\ref{h-ODE2}) is designated by green.

\end{document}